\documentclass[11pt]{amsart}
\usepackage{geometry}                
\usepackage{graphicx}
\usepackage{amssymb}
\usepackage{mathabx}
\usepackage{epstopdf}
\usepackage{amsthm} 
\usepackage{enumerate}
\usepackage[all]{xy}
\usepackage[linktocpage]{hyperref}

\newtheorem{thm}{Theorem}

\numberwithin{thm}{section}
\newtheorem{lemma}[thm]{Lemma}

\newtheorem{proposition}[thm]{Proposition}

\newtheorem{corollary}[thm]{Corollary}

\theoremstyle{remark}
\newtheorem{rmk}{Remark}[thm]

\theoremstyle{definition}
\newtheorem{definition}[thm]{Definition}

\newcommand{\QED}{\hfill $\Box$}
\newcommand{\norm}[1]{\left\|#1\right\|}
\newcommand{\abs}[1]{\left\vert#1\right\vert}
\newcommand{\set}[1]{\left\{#1\right\}}
\newcommand{\overl}[1]{\overline{#1}}

\newcommand{\mbb}{\mathbb}

\newcommand{\R}{\mathbb{R}}

\newcommand{\C}{\mathbb{C}}

\newcommand{\Z}{\mathbb{Z}}
\newcommand{\Q}{\mathbb{Q}}

\newcommand{\utp}{\widetilde{type}}
\newcommand{\mlaur}{\Z\left[z_1,z_1^{-1}\right]\oplus\cdots\oplus\Z\left[z_k,z_k^{-1}\right]}

\newcommand{\Tl}[1]{Tail_{h}\left(#1\right)}
\newcommand{\Ct}[1]{Center_{h}\left(#1\right)}
\newcommand{\Hd}[1]{Head_{h}\left(#1\right)}
\newcommand{\lp}[2]{\overl{#1}=(#1_1,\ldots,#1_{#2})}
\newcommand{\dv}{\textbf{div}}
\newcommand{\Dv}{\textbf{Div}}

\title{Rational Growth and Almost Convexity of Higher-dimensional Torus Bundles}
\author{Corey Bregman}
\begin{document}
\maketitle

\begin{abstract} Given a matrix $A\in SL(N,\Z)$, form the semidirect product $G=\Z^N\rtimes_A \Z$ where the $\Z$ factor acts on $\Z^N$ by $A$.  Such a $G$ arises naturally as the fundamental group of an $N$-dimensional torus bundle which fibers over the circle.  In this paper we prove that if $A$ has distinct eigenvalues not lying on the unit circle, then there exists a finite index subgroup $H\leq G$ possessing rational growth series for some generating set. In contrast, we show that if $A$ has at least one eigenvalue not lying on the unit circle, then $G$ is not almost convex for any generating set.   

\end{abstract}

\section{Introduction}\label{intro}

Let $G$ be a finitely generated group with generating set $S=S^{-1}$.  We equip $G$ with the word-norm $\|\cdot\|_S$ and word-metric $d_S(\cdot,\cdot)$ coming from $S$.   
This metric is one of basic objects of modern geometric group theory, and it is well-known that many important coarse geometric invariants of $G$ can be defined in a way which is independent of the generating set.  
Not as well understood are metric properties of $G$ which \emph{depend} on a particular generating set.  In this paper, we will explore two such properties, rational growth and almost convexity, for fundamental groups of higher-dimensional torus bundles.

\subsection{Rational Growth and Almost Convexity}
 Let $b_{G,S}(n)$ denote the number of group elements in the ball of radius $n\geq0$, centered on the identity in $G$.  One way to understand the growth of $G$ is to consider the sequence $\{b_{G,S}(n)\}$. One defines the growth series of $G$ with generating set $S$ to be the formal power series:
\[P_{G,S}(z) = \sum_{n=0}^{\infty} b_{G,S}(n)z^n.\]
This series always has a positive radius of convergence, and one can try to find a generating function associated to it.  We say that a pair $(G,S)$ has \emph{rational} (resp. \emph{algebraic} or \emph{transcendental}) \emph{growth} if the above series has a generating function which is rational (resp. algebraic or transcendental) over the field $\Q(z)$. We remark that having rational growth series is equivalent to having a finite system of linear recursion equations among the $\{b_{G,S}(n)\}$.

A related notion is \emph{almost convexity}, first defined by Cannon in \cite{CAC87}. A pair $(G,S)$ is said to be \emph{almost convex} if for every $k\geq 0$ there exists a constant $N(k)$ so that the following holds. If $g,h\in G$ satisfy $\|g\|_S=\|h\|_S \mbox{ and } d_S(g,h)=k$, then there is a path in the ball of radius $n$ between $g$ and $h$ whose length is at most $N(k)$.  

Cannon proved that when a pair $(G,S)$ is almost convex the group has solvable word problem, by showing that the Cayley graph can be constructed algorithmically.  More precisely, the ball of radius $n$ in the Cayley graph can be constructed from the ball of radius $n-1$ by applying a finite number of rules.  One might expect some relationship between almost convexity and rational growth; to wit, if one can build the ball of radius $n$ from balls of smaller radius, the size of the ball of radius $n$ ought to depend on the sizes of smaller balls in some controlled way and vice versa.

The groups we will consider are polycyclic extensions of $\Z$ by $\Z^N$ formed in the following way. Let $A \in SL(N,\Z)$ be any matrix. Form the semidirect product $G\cong \Z^N\rtimes_A\Z$ where the generator of the $\Z$ factor conjugates $\Z^N$ to $A\cdot \Z^N$.  Such groups arise naturally as fundamental groups of higher-dimensional torus bundles which fiber over the circle.  In connection with rational growth, we prove:

\begin{thm}\label{main1}
 Let $A\in SL(N,\Z)$ be a matrix with distinct eigenvalues not lying on the unit circle, and let $G\cong \Z^N\rtimes_{A}\Z$ be the extension of $\Z^N$ where the $\Z$ factor acts on $\Z^N$ by $A$.  Then there exists a finite index subgroup of $G$ with finite generating set $S$ whose growth series is rational.

\end{thm}

Theorem \ref{main1} generalizes the work of Parry \cite{Pa07} and Putman \cite{Pu06}, who establish similar results when $N=2$.  In contrast, we also prove the following, extending the main theorem of Cannon, Floyd, Grayson and Thurston \cite{CaFl89} to higher dimensions 

 \begin{thm} \label{main2} Let $A\in SL(N,\Z)$ be a matrix with at least one eigenvalue not lying on the unit circle and let $G\cong \Z^N\rtimes_{A}\Z$ be the extension of $\Z^N$ where the $\Z$ factor acts on $\Z^N$ by $A$.  Then $G$ is not almost convex with respect to any generating set.  

\end{thm}

\begin{rmk} An unpublished 2008 preprint \cite{AW08} of Andrew Warshall claims to prove the a similar result to Theorem \ref{main2}. However, Warshall requires the automorphism induced by $A$ to be hyperbolic, \emph{ i.e.} for all eigenvalues of $A$ to lie off the unit circle, whereas our theorem only requires at least one.  Moreover, the method of proof in Warshall's result is entirely different from our own. 


\end{rmk}
%

\subsection{Background}
In general, it is a difficult and subtle combinatorial problem to compute the growth series of some pair $(G,S)$. There are, however, some general obstructions to a group possessing rational growth.  For example, it is known that the coefficients of a Taylor series expansion of a rational function grow either polynomially or exponentially.  Hence a group of intermediate growth, such as the Grigorchuk group \cite{Gri83}, cannot possess a rational growth series for any generating set.  Moreover, having a rational growth series implies solvable word problem, so any finitely generated group with unsolvable word problem cannot have rational growth for any generating set \cite{Ca80}.  These obstructions, though fairly general, do not seem to apply to the types of groups one usually encounters. 

Although rational growth is generating set dependent \cite{St96}, it is known that for certain large classes of groups the growth series is rational with respect to some  generating set.  For example, an exercise in Bourbaki outlines a proof that a finitely generated Coxeter group has rational growth series with respect to a Coxeter system \cite{Bo68}.  Cannon proved that any cocompact hyperbolic group has rational growth with respect to \emph{any} generating set; his argument easily generalizes to any word-hyperbolic group \cite{Ca84}.  In fact, in \cite{NS95} Neumann and Shapiro prove that, among other things, for any geometrically finite hyperbolic group, any generating set can be extended to a generating set possessing rational growth series.   

For finitely generated abelian groups, rational growth with respect to any generating set follows from a classical result of Hilbert.  Benson has proven that for any virtually abelian group, rational growth is also independent of generating set \cite{Be83}.  In light of Bieberbach's theorem, even in the last case existence of rational growth for some generating set seems to be strongly tied to geometric properties of the group.  Similar results hold for almost convexity as well; namely, Thiel has shown that almost convexity depends on the generating set \cite{CT91}, but it is known that virtually abelian and word-hyperbolic groups are almost convex for every generating set \cite{CAC87}.

Following the work of Perelman and Thurston, another natural class of groups to consider is that of fundamental groups of compact geometric 3-manifolds (See \cite{Th97} and \cite{Sc83}). Existence of rational growth series and almost convex generating sets for cocompact subgroups of the  geometries $\mathbb{E}^3$, $\mathbb{S}^3$, $\mathbb{H}^3$, $\mathbb{H}^2\times\R$, and $\mathbb{S}^2\times \R$ are covered by the above theorems of Benson and Cannon.  For the remaining three geometries, \textbf{Nil}, $\widehat{PSL}_2(\R)$, and \textbf{Sol}, some partial results on rational growth are known (for \textbf{Nil} see \cite{Be84}, \cite{Sh89},  \cite{St96},  and the recent paper of Duchin and Shapiro \cite{DS14}, for $\widehat{PSL}_2(\R)$ see \cite{Sh94}).  On the other hand, for every 3-manifold geometry except \textbf{Sol}, almost convexity holds for some generating set \cite{SSt95}, while  \textbf{Sol} groups are not almost convex for any generating set  \cite{CaFl89}.  


Compact manifolds admitting \textbf{Sol}-geometry are precisely torus bundles over the circle with Anosov monodromy.  These groups therefore have presentations as semidirect products $\Z^2\rtimes\Z$, where the generator of the $\Z$-factor acts on $\Z^2$ via some matrix in $A\in SL(2,\Z)$ whose trace is greater than 2.  When the trace is even, Parry has computed some examples of the growth series and shown it to be rational \cite{Pa07}.  Putman, on the other hand, has shown that every such group has a finite index subgroup whose growth series is rational \cite{Pu06}.  The results of Parry and Putman are striking in light of the fact that \textbf{Sol} groups are not almost convex.


The three-dimensional \textbf{Sol} examples fit into a larger class of torsion-free abelian-by-cyclic groups.  By a result of Bieri and Strebel, such groups are exactly of the form $\Z^N\rtimes_{A}\Z$, where $A$ is a non-singular matrix with $\Z$-coefficients \cite{BS78}. These groups divide naturally into polycyclic and non-polycyclic types, and are polycyclic only when $A\in GL(N,\Z)$.  In the non-polycyclic case, Brazil \cite{Bra94} and Collins-Edjvet-Gill \cite{Edjvet94} have demonstrated rational growth for the solvable Baumslag-Solitar groups $BS(1,n)$ with the usual generating set. In this context, the groups we consider in this paper form a large subset of all torsion-free polycyclic abelian-by-cyclic groups which contains all of the 3-dimensional \textbf{Sol} examples.  

%
%
%
%

\subsection{Outline}
The paper is divided into two parts, corresponding to the two theorems above.  The structure of each part is as follows.

\textbf{Part I}: In Section \ref{tbundles},  we describe the construction of higher dimensional torus bundles over the circle, and prove that we can always pass to a finite index subgroup so that the corresponding torus bundle has a particularly nice presentation.  In Section \ref{RegLang}, we review some relevant definitions and theorems from the theory of regular languages, and recall the falsification-by-fellow-traveler property (FFTP), first introduced in \cite{Pu06}, and inspired by a similar idea in \cite{NS95}.  The FFTP will allow us to deduce that certain partitions of regular languages have rational growth. We introduce the notion of types and heights in Section \ref{Lang}, and then describe a family of regular languages $\{\Lambda_n\}$, each of which surject onto $G$.  Finally, in Section \ref{FFTP}, we prove that for $n$ sufficiently large, $\Lambda_n$ has the FFTP, and deduce Theorem \ref{main1}.  

\textbf{Part II}: In Section \ref{AC}, we review the definition of almost convexity.  In Section \ref{Solv} we revisit the torus bundles introduced in Section \ref{tbundles} by representing them as certain geometric quotients of two-step solvable Lie groups. Embedding the groups above as lattices in Lie groups allows us to exploit the asymptotic geometry of the Lie groups to prove Theorem \ref{main2}.  

\textbf{Acknowledgements.} I would like to thank my advisor, Andrew Putman, for introducing me to this problem and for the many helpful discussions we had while the paper was being written.  We would also like to thank Spencer Dowdall for pointing out to us the result in \cite{AW08}.

\begin{center}\Large \textbf{Part I: Rational Growth}
\end{center}
\section{Higher-dimensional torus bundles}\label{tbundles}
Let $N>1$ be an integer.  We consider $N$-dimensional torus bundles over the circle defined in the following way.  Let $\mbb{T}^N$ be the standard square $N$-dimensional torus.  Up to homotopy, an orientation-preserving homeomorphism of $\mbb{T}^N$ can be represented by a matrix  $A \in SL(N,\mbb{Z})$.  The linear map $T_A:\R^N\rightarrow\R^N$ where $T_A(x)=Ax$ takes the lattice to the lattice, and hence defines a map on the quotient $\overline{T}_A:\mbb{T}^N\rightarrow \mbb{T}^N$. We can take $\overline{T}_A$ as a representative in each homotopy class.  

Each $A \in SL(N,\mbb{Z})$ also defines a torus bundle over the circle by taking $\mbb{T}^N\times I$ and identifying $\mbb{T}^N\times\{0\}$ with  $\mbb{T}^N\times\{1\}$ via the homeomorphism $\overline{T}_A$.  That is, for each $A \in SL(N,\mbb{Z})$ we define the bundle $X_A$ as the quotient space
\[X_A :=\mbb{T}^N\times I/\{(x,0)\sim(\overline{T}_A(x),1)\}.\]

The fundamental group $\pi_1(X_A,*)$ is an HNN-extension of $\Z^N$ by $\Z$.  By the Seifert-van Kampen theorem \[\pi_1(X_A,*)= \langle \Z^N,t |t  \Z^Nt^{-1}=A\Z^N\rangle \cong\Z^N\rtimes_A\Z.\]  

An easy consequence of Gauss' lemma is that if a matrix $A \in SL(N,\mbb{Z})$ has an integral eigenvalue $\lambda$, then $\lambda=\pm 1$.  It follows that if $\lambda$ is a root of an integer, then $\lambda$ is a root of unity.   Our results will apply to those $A$ with distinct eigenvalues not lying on the unit circle.  Thus the matrices we consider will not have any integral eigenvectors.  

\subsection{Block Rational Canonical Form}
Let $A\in SL(N,\Z)$ be a matrix with distinct eigenvalues not lying on the unit circle, and let $\{e_1,\ldots,e_N\}$ be the standard basis for $\Q^N$. The characteristic polynomial for $A$, denoted $char(A)$, factors into distinct irreducible polynomials over $\mbb{Q}$, none of which are linear:\[char(A)=p_1\cdots p_k.\] Set $d_i=\deg(p_i)$, so that $\sum_i d_i=N$.  Since the eigenvalues of $A$ are assumed distinct, all the $p_i$ are distinct as well.  As a $\Q[t]$-module, $V=\Q^N$ decomposes as a direct sum of cyclic modules:\[V\cong \Q[t]/(p_1(t))\oplus\cdots\oplus\Q[t]/(p_k(t)).\] By a change of basis, we can put our matrix in a form which mirrors the decomposition above:
\[A=\left(\begin{array}{cccc}
					A_1 &0&\cdots&0\\
					0& A_2&\cdots&0\\
					0&\cdots&\ddots&0\\
					0&\cdots&0&A_k	

				\end{array}\right)\] 
where the $A_i$ are $(d_i\times d_i)$-blocks and $char(A_i)=p_i$. Now writing $V\cong V_1\oplus \cdots\oplus V_k$, we see that $A$ acts by $A_i$ on $V_i$.  Choose an integral vector $v_i\in V_i$. Since the $p_i$ are irreducible over $\Q$, $\{v_i, A_iv_i,\ldots,A_i^{d_i-1}v_i\}$ are linearly independent over $\Q$ and $p_i(A_i)v_i=0$.  By changing our basis to $\{v_1,\ldots,A_1^{d_1-1}v_1,\ldots, v_k,\ldots,A_k^{d_k-1}v_k\}$, we can assume each $A_i$ block is in rational canonical form. 

\begin{definition} A matrix $A$ in form described above is said to be in \textbf{block rational canonical form}.
\end{definition}
\subsection{Finite index subgroups of $G$}
Although these changes of basis in the previous section occur over $\Q^N$, they can be achieved over $\Z^N$ by passing to a finite index subgroup.  Form the HNN-extension $G\cong \Z^N\rtimes_A\Z$ as above. Let $\{e_1,\ldots,e_N\}$ be the standard basis for $\Z^N$, and $\{t\}$ the generator of $\Z$ in the HNN-extension. We may pass to a finite index subgroup of $G$ by alternately:
\begin{enumerate}
	\item Replacing $\Z^N$ by a finite index subgroup of $\Z^N$.  This amounts to choosing a different basis for $\Q^N$.
	\item Replacing $t$ by $t^k$ for some $k>1$.  This is equivalent to replacing $A$ by $A^k$.  
\end{enumerate}

In order to prove rational growth of a subgroup of $G$, we will need each of the factors $p_i$ of $char(A)$ to satisfy a special property.   First write $p_i(z) = \sum_ja_{ij}z^j$. Each $p_i$ must satisfy:
  
$(\ast)$ There exists $j_0$ so that $\displaystyle |a_{ij_0}|>\sum_{j\neq j_0}|a_{ij}|$.

It will be convenient to define two norms on the set of Laurent polynomials $\Z[z,z^{-1}]$.  Let $\displaystyle f(z) = \sum_{i=-n_1}^{n_2}a_iz^i$.  Then we define \[|f|_{\infty} = \max_i\{|a_i|\}\] and \[ |f|_{1} = \sum_i|a_i|.\] Now $(\ast)$ can be expressed concisely as:

$(\ast)$ For each $i$, $|p_i|_{\infty}>\frac{1}{2}|p_i|_1$.

The main result of this section is the following lemma:
\begin{lemma}\label{finind} Let $G$ be an HNN extension as above.  By passing to a finite index subgroup of $G$ we may assume that \begin{enumerate}
	\item $A$ is in block rational canonical form over $\Q$ with coefficients in $\Z$.
	\item Each irreducible factor $p_i$ of the $char(A)$ possesses a maximal coefficient $M_i$ satisfying $|M_i|=|p_i|_{\infty}>\frac{1}{2}|p_i|_{1}$.  
	\end{enumerate}

\end{lemma}

\footnotesize{PROOF}: 
\normalsize The proof of (1) is nearly complete from the discussion above.  Observe that the change of basis matrix may be chosen with entries in $\Z$.  This is easily achieved by clearing denominators in the change of basis matrix over $\Q$. Finally, note that because $char(A)$ has all distinct roots, it is possible to find infinitely many powers $A^k$ such that $char(A^k)$ has distinct roots. This guarantees that after taking powers of $A$ we can still obtain a block rational canonical form.  

For (2), first assume that $char(A^n)$ is irreducible for all powers $n$.  Order the roots of $\lambda_1,\ldots,\lambda_N$ of $A$ according to their modulus: \[|\lambda_1|\geq\cdots \geq |\lambda_r|>1>|\lambda_{r+1}|\geq \cdots \geq |\lambda_N|.\]  
Here we make use of the fact that no root lies on the unit circle and that the product of the roots is 1.  Set  $l_i=|\lambda_1\cdots\lambda_i|$.  Then \[l_1<l_2<\cdots <l_r \mbox{ and }l_r>l_{r+1}>\cdots >l_N.\] Since $char(A^n) = \prod_i (z-\lambda_i^n)$, we see that the $k$th coefficient $c_{n,k}$ of $char(A^n)$ can be written:  
\[c_{n,k}=\sum_{I = (i_1,\ldots, i_k)}\lambda_I^n=(\lambda_1\cdots \lambda_k)^n \left(1+\sum_{I\neq(1,\ldots,k)}\frac{\lambda_I^n}{(\lambda_{1} \cdots \lambda_{k})^n}\right)\] where $1\leq i_1<\cdots<i_k\leq N$ and $\lambda_I=\lambda_{i_1}\cdots\lambda_{i_k}$. Then\[|c_{n,k}| = l_k^n\cdot \left|1+\sum_{I\neq(1,\ldots,k)}\frac{\lambda_I^n}{(\lambda_{1} \cdots \lambda_{k})^n}\right|.\]  
By choosing a sufficiently high power of $A$, the sum on the right can be made arbitrarily close to 1. It follows that  $|c_{k,n}|$ is asymptotic to $l_k^n$.  Since $l_r>l_k$ for all $k\neq r$, for $n$ sufficiently large $c_{r,n}$ can be made large enough to satisfy the conclusion of (2).  

Taking into account the fact that $char(A^n)$ may not be irreducible, consider all possible combinations of roots which multiply to give $\pm1$.  For each of these, find an integer to satisfy the conclusion of (2).   Choose $n_{\text{max}}$ to be the product of all of these integers.  This choice of $n$ will guarantee that even if $char(A^n)$ factors, each irreducible factor must have a coefficient which satisfies condition (2) of the lemma.         \QED

\section{Weighted Sets and Regular Languages}\label{RegLang}

\subsection{Weighted Sets} We introduce the formal notions of \emph{weightings} and \emph{weighted partitions} on arbitrary sets to streamline our constructions later on.

Let $X$ be a set.  A \emph{weighting} on $X$ is a function $\norm{\cdot}:X\rightarrow\Z_{\geq0}$.  A partition $P$ of $X$ is a decomposition of $X$ into pairwise disjoint subsets: \[X=\coprod_{A\in P}A.\]
Denote by $X/P$ the set of equivalence classes of $X$ modulo $P$.  For every $x\in X$ we denote by $\overl{x}$ the equivalence class of $x$ determined by $P$.  There is a natural surjective map $X\rightarrow X/P$, sending $x\mapsto \overl{x}$. If $X$ has a weighting then the quotient $X/P$ inherits a weighting via \[\norm{\overl{x}} = \min\{\norm{y}\vert y\in X, \overl{y} =\overl{x}\}.\]
We call $X/P$ a \emph{weighted partition}.

\begin{definition} Let $(X_1,\norm{\cdot}_1)$, $(X_2,\norm{\cdot}_2)$ be weighted sets.  A bijection $\psi:X_1\rightarrow X_2$ is called a \textbf{near-isometry} if there exists a constant $c$ such that  for each $x \in X_1$,  $\norm{\psi(x)}_2 =\norm{x}_1+c$.
\end{definition}

\subsection{Regular languages} In this section we recall results from the theory of regular languages that will be relevant for the proof of Theorem \ref{main1}.  For more details, consult \cite{ECHLPT92}, Ch. 1.  We designate a finite set $A$ as an \emph{alphabet} of symbols.  Each element of $A$ is called a \emph{letter}.  We can form \emph{words} or strings over $A$ by concatenating letters. We denote by $\epsilon$ the empty string and by $A^{*}$ the set of all strings. The length $|w|$ of a word $w\in A^{*}$ is the number of letters in $w$.  
\begin{definition} A \textbf{language} over $A$ is a subset $L\subset A^{*}$.
\end{definition}  

Certain languages behave better than others.  We will be interested in finding regular languages that are in bijection with elements in the group $G$.  We can think of the language as a `normal form' for group elements.  Regular languages are those which can be recognized by a finite-state automaton:

\begin{definition} A \textbf{finite-state automaton} (FSA) is a 5-tuple $M=(S,A,\mu, F, s_0)$ where
\begin{itemize}
		\item $S$ is a finite set of states.
		\item $A$ is an alphabet.
		\item $\mu: S\times A\rightarrow S$ is the \emph{transition} function.
		\item $F\subset S$ is the collection of final or \emph{accept} states.
		\item $s_0\in S$ is the initial state.  
\end{itemize}

\end{definition}

We think of $M$ as a machine which can read strings one letter at a time.  If $M$ reads the symbol $a\in A$ while in the state $s\in S$, it proceeds to the state $\mu(s,a)$ and reads the next letter.  We say $M$ \emph{accepts} a string if $M$ transitions to one of the accept states upon reading the final letter in the string. Denote by $L_M$ the language of all strings accepted by $M$. 

\begin{definition} A language $L$ over $A$ is said to be \textbf{regular} if $L=L_M$ for some FSA $M$.  
\end{definition}

In the course of the proof, it will be important to encode first-order predicate relations between strings in a regular language $L$, e.g. $|w_1|<|w_2|$.  The best way to do this is to introduce the notion of an $n$-variable language.

\begin{definition}Let $A$ be an alphabet not containing the symbol $\{\$\}$ and let $L\subset A^{*}$ be a language. We form a new alphabet $A^{\$}= A\cup \{\$\}$ and define a language $L_n^{\$}\subset \prod_{i=1}^n(A^{\$})^*$ as follows.  Consider any $n$-tuple $(w_1,\ldots,w_n) \in \prod_{i=1}^nL$.  Let $m$ be the length of the longest word in this $n$-tuple.  If some $w_i$ does not have length $m$, we pad $w_i$ by adjoining copies of $\$$ on the end of it until it has length $m$. We thus obtain an $n$-tuple $(w_1,\ldots,w_n)^{\$}\in \prod_{i=1}^n(A^{\$})^*$ and define \[L_n^{\$} = \set{(w_1,\ldots,w_n)^{\$} \vert (w_1,\ldots,w_n) \in \prod_{i=1}^nL}.\]

$L_n^{\$}$ is called the (padded) \textbf{$n$-variable language} over $L$ (with padding $\$$).
\end{definition}  

\begin{rmk} $L_n^{\$}$ is a language when considered over the alphabet $\prod_{i=1}^n(A^{\$})^*$.  Note that if $L\subset A^{*}$ is regular, then $L_n^{\$}$ is regular over $\prod_{i=1}^n(A^{\$})^*$.
\end{rmk} 

With this definition it is possible to form first-order predicate relations among regular languages. We have the following useful result:

\begin{thm}\label{ECHprop}$($\cite{ECHLPT92}, Proposition $1.1.4$, Theorem $1.2.8$, Corollary $1.4.7$$)$
The set of regular languages is closed under:
\begin{enumerate}
	\item First order predicates: complementation, union, intersection, $\exists$, $\forall$, concatenation, and Kleene star. 
	\item Reversal:  if $L$ is a regular language, then the language \[\text{rev}(L)= \{w^r=a_n\cdots a_1\vert w=a_1\cdots a_n \in L\}\] is regular.
\end{enumerate}

\end{thm}

We can put a weighting on $A^{*}$ by starting with a function $\phi:A\rightarrow\Z_{\geq0}$ and extending it to a weighting on words $w=a_1\cdots a_n\in  A^{*}$ via $\norm{w}=\sum_{i=1}^n\phi(a_i)$.   Given a weighting $\phi$ and a language $L\subset A^{*}$ we define the \emph{growth series} of $L$ with weighting $\phi$ to be the formal power series \[G_{\phi}(L) = \sum_{i=0}^{\infty}c_iz^i\] where \[c_i=\#\{w\in L\vert \phi(w)=i\}.\]

\begin{rmk} If $\phi$ is identically 1, the weighting is just the ordinary word length as defined above.
\end{rmk}

In some cases the growth series of a language may represent the Taylor series expansion about zero of a function which is rational, algebraic, or transcendental over $\Q\left(z\right)$.  The relatively simple structure of the finite state automaton ensures that the growth series for a regular language is rational with respect to any weighting.  
\begin{thm}\label{regrat}  Let $L$ be a regular language and $\phi:A\rightarrow\Z_{\geq0}$ a weighting. Then $G_{\phi}(L)$ is a rational function.
\end{thm}

\subsection{The Falsification-by-fellow-traveler Property}  Let $L$ be a regular language with weighting $\phi$.  Suppose $P$ is a partition for $L$ and let $\pi:L\rightarrow L/P$ be the quotient map.  Following Putman we define what it means for the pair $(L,P)$ to have the \emph{falsification-by-fellow-traveler} property.  For proofs of the results in this section, see \cite{Pu06}.

\begin{definition} We say that $L/P$ has a \textbf{regular cross-section} if there is a regular sublanguage $L'\subset L$ and a surjective map $\sigma:L/P\rightarrow L'$ so that $\pi\circ\sigma=\text{Id}_{L/P}$. $L'$ is called \textbf{minimal} if it satisfies\[ \norm{\sigma(A)} = \min\{\norm{x}\vert x\in A\}\]
for every $A\in P$.
\end{definition}

\begin{definition} A regular language $R\subset L\times L$ is an \textbf{acceptor} for a partition $P$ of $L$ if \[(w,w')\in R\Rightarrow \overl{w}=\overl{w}' \text{ and } (w',w)\in R.\]
\end{definition}

With these two definitions in mind, we can state the falsification-by-fellow-traveler property:
\begin{definition} The pair $(L,P)$ with acceptor $R$ has the \textbf{falsification-by-fellow-traveler property (FFTP)} if there is a constant $K>0$ and a regular sublanguage $L'\subset L$ which contains at least one minimal size representative of each equivalence class in $P$ and which satisfies the following.  If $w\in L'$ is not a minimal representative modulo $P$, then there exists $w'\in L$ with:
\begin{itemize}
	\item $(w,w')\in R$
	\item $\norm{w'}<\norm{w}$
	\item For any $i$, let $s$ and $s'$ be the length $i$ initial segments of $w$ and $w'$, respectively.  Then $|\norm{s}-\norm{s'}|\leq K$. In this case, $w$ and $w'$ are said to $K$-\emph{fellow-travel}.
\end{itemize}
We also require that if $w,w'\in L'$ are two minimal size representatives for the same equivalence class then $(w,w')\in R$.  
\end{definition}
The usefulness of this property is the following theorem, proved in \cite{Pu06}:
\begin{thm}\label{PutFFT} $($\cite{Pu06}, Theorem $3.1)$ Let $L$ be a weighted regular language, $P$ a partition on $L$. Suppose the pair $(L,P)$ with acceptor $R$ has the FFTP.  Then $P$ has a regular minimal cross-section. In particular, $L/P$ has rational growth.
\end{thm}

\section{The Language}\label{Lang}
\subsection{Types and Heights}  Thinking of $G$ as the fundamental group of a fiber bundle, we will decompose group elements of $G$ into their fiber and base components.  The fiber component will be called the \emph{type} and the base component will be called the \emph{height}.  This construction goes back to work of Grayson in his thesis \cite{Gray83}, and was adopted by subsequent papers on rational growth in cocompact subgroups of $Sol$ (\cite{Pa07},\cite{Pu06}).  Using the same notation as in Section \ref{tbundles} we denote by $A \in SL(N,\mbb{Z})$ a matrix in block rational canonical form with distinct eigenvalues not lying on the unit circle.  Let $char(A)=p_1\cdots p_k$ be the factorization of the characteristic polynomial of $A$ into irreducibles.  We assume that the conclusion of Lemma \ref{finind} holds for $A$.

Let $\{a_1,\ldots,a_k\}$ be a set of generators for $\Z^k$.  Since $G$ is an HNN-extension, we have the following commutative diagram:
\[  \xymatrix{
1  \ar[r] &  K\ar[d] \ar[r]^{\iota'}        &    \Z^k*F_{\{t\}} \ar[d]^\rho \ar[r]^-{\pi'}&   F_{\{t\}}\ar[d]^-\cong \ar[r] & 1  \\
1\ar[r]&  \Z^N  \ar[r]^{\iota}   &  G\ar[r]^-\pi  &   \Z \ar[r]& 1
  }
\]
where $F_{\{t\}}$ is the free group on the generator $\{ t\}$, and $K$ is the kernel of the projection map onto the factor $F_{\{t\}}$.  The kernel $K$ can easily be identified with the free product $\displaystyle *_{i\in \Z}\Z^k$ with generating set \[\{t^ia_jt^{-i}|1\leq j \leq k, i\geq0\}.\]  

To see this, realize $\Z^k*\Z$ as the fundamental group of $\mbb{T}^k\vee\mbb{T}^1$, and observe that $N$ corresponds to the cover homeomorphic to a copy of $\R$ with a $k$-torus $\mbb{T}^k$ wedged on at every integer point.  Note that since $\Z^N$ is abelian, the map $K\rightarrow\Z^N$ factors through the abelianization of $K$, which we identify with $\mlaur$ by sending $t^ia_jt^{-i}\mapsto z_j^i$.  

Because $F_{\{t\}}$ is free, the top sequence splits, and we can therefore represent $\Z^k*F_{\{t\}}$ as the semi-direct product $K\rtimes F_{\{t\}}$.  This allows us to represent each element of $\Z^k*F_{\{t\}}$ uniquely as a pair $(v,t^h)$, where $v$ is an element of $K$ and $t^h$ is an element of $F_{\{t\}}$.

Now let $w$ be a freely reduced word in the generators $\{a_1,\ldots,a_k,t\}$, and denote by $\overl{w}$ the group element in $\Z^k*F_{\{t\}}$ it represents.  As an element of the above semi-direct product, we can write $\overl{w}=(v,t^h)$. 

\begin{definition} The \textbf{height} of $w$ is $ht(w)=h=\pi\circ\rho(\overl{w})\in \Z$, and the \textbf{type} of $w$ is $type(w)=\iota^{-1}\circ \rho(v)\in \Z^N$.  The \textbf{unreduced type} of $w$, denoted $\widetilde{type}(w)$, is the image of $v$ in the abelianization $K^{ab}$.  
\end{definition}
\begin{rmk}
The manifold $X_A$ being the total space of a fiber bundle over the circle $S^1$, we think of $ht(w)$ somewhat paradoxically as the horizontal or base component of $w$, and $type(w)$ as the vertical or fiber component.  
\end{rmk}

The unreduced type $\utp(w)\in \mlaur$ is a tuple of Laurent polynomials, which will frequently be denoted by $\overl{t}=\left(t_1(z_1),\ldots,t_k(z_k)\right)$.  The norms $\abs{\cdot}_{\infty}$ and $\abs{\cdot}_{1}$ extend naturally to tuples of Laurent polynomials.  Following Putman, we define functions $\Tl{f}$, $\Ct{f}$, and $\Hd{f}$ for a Laurent polynomial $f(z)$ which depend on $h$ as follows:

For a Laurent polynomial \[f(z)= \sum_ic_iz^i\in \Z\left[z,z^{-1}\right]\] and $h \in \Z$ an integer we partition $f$ into three pieces depending on $h$.  If $h\geq 0$, \[\Tl{f}:=\sum_{i=-\infty}^{-1}c_iz^i,\]

\[\Ct{f}: = \sum_{i=0}^{h}c_iz^i,\]

\[\Hd{f}:=\sum_{i=h+1}^{\infty}c_iz^i,\]

For $h\leq0$, we extend this definition symmetrically via \[\Tl{f}(z^{-1}) = Head_{-h}(f(z^{-1})),\]
\[\Ct{f}(z^{-1}):=Center_{-h}(f(z^{-1})),\] 
\[\Hd{f}(z^{-1}) := Tail_{-h}(f(z^{-1})). \]

We also define corresponding lengths $T_h(f)$, $H_h(f)$ as \[T_h(f) :=\max\{\abs{i}\vert i\leq0,c_i\neq0\},\] \[H_h(f) :=\max\{i-h\vert i\geq h,c_i\neq0\}\] if $h\geq 0$ and as  \[T_h(f) :=\max\{\abs{i}-\abs{h}\vert i\leq h,c_i\neq0\},\] \[H_h(f) :=\max\{i\vert i\geq 0,c_i\neq0\}\] if $h\leq0$. 

We extend each of these functions to tuples of Laurent polynomials $\overl{t}=\left((t_1(z_1),\ldots,t_k(z_k)\right)$ via:
\begin{enumerate}[(a)]

\item $\Tl{\overl{t}} := \left(\Tl{t_1},\ldots,\Tl{t_k}\right)$ and similarly for $\Ct{\overl{t}}$ and $\Hd{\overl{t}}$.
\item $T_h(\overl{t}):=\max\left\{T_h(t_1),\ldots,T_h(t_k)\right\}$ and similarly for $H_h(\overl{t})$.
\end{enumerate}

We now determine the shortest length of a word $w$ with $\utp(w) = \overl{t}$ and $ht(w)=h$:
\begin{proposition}\label{lengthprop}
Let $h$ be a height and $\overl{t}=(t_1(z_1),\ldots, t_k(z_k))$ be an unreduced type.  Then any shortest word with this unreduced type and height has length
		\[2T_h(\overl{t})+2H_h(\overl{t})+|h| +|\overl{t}|_1.\]
\end{proposition}
\footnotesize{PROOF}: 
\normalsize  Note that the $z_j$ all commute. The proof is the same as in \cite{Pu06}, Theorem 4.2.  \QED

If $w$ is a freely reduced word in the generators, the set of pairs $(\utp(w),ht(w))$ certainly surjects onto the set of elements of $G$, but is not in bijection with it.  Even if we restrict the set of $w$ to those with minimal length unreduced types, there will still be more than one representative for each element of $G$.  

In order to get a handle on the discrepancy exactly, we must introduce a variant of polynomial division for tuples of Laurent polynomials.  We say that $\lp{f}{k}$ divides $\lp{g}{k}$, denoted $\overl{f}|\overl{g}$, if $f_i|g_i$ for all $1\leq i \leq k$.  Denote by $\overl{p}=(p_1,\ldots,p_k)$ the tuple of polynomials whose components are the irreducible factors of $char(A)$. We have  

\begin{proposition}\label{grpequal}

Let $w_1$, $w_2$ be words in the generators for $G$.  Then  $w_1=_G w_2$ if and only if $ht(w_1)=ht(w_2)$ and $\overl{p}=(p_1,\ldots,p_k)$ divides $\utp(w_1)-\utp(w_2)$.  

\end{proposition}
\footnotesize{PROOF}: 
\normalsize Recall that our chosen module generating set in block rational canonical form consists of vectors $\{v_1,\ldots, v_k\}$. First observe that if $\utp(w) = \lp{t}{k}$ then \[type(w)= t_1(A)\cdot v_1+\cdots +t_k(A)\cdot v_k\] by our choice of basis for block rational canonical form.

Since the minimal polynomial for $A$ divides the characteristic polynomial, the two conditions are clearly sufficient.  On the other hand, note that if $w_1=_Gw_2$ then we must have $ht(w_1)=ht(w_2)$ by projecting onto $\Z$. Moreover, the fact that the $p_i$ are distinct and irreducible implies that the characteristic polynomial is the minimal polynomial hence the two conditions are necessary.\QED

Define a partition $P$ on $S = \mbb{Z}\left[z_1,z_1^{-1}\right] \times\cdots\times\mbb{Z}\left[z_k,z_k^{-1}\right] \times \Z$ via \[(\overl{t}_1,h_1)\sim(\overl{t}_2,h_2) \text{ iff } h_1=h_2 \text{ and } (p_1,\ldots,p_k) \text{ divides }\overl{t}_1-\overl{t}_2.\]

Define a weighting on $S$ by $\norm{(\overl{t},h)} = 2T_h(\overl{t})+2H_h(\overl{t})+|h| +|\overl{t}|_1 +1$.
Then the previous proposition implies:
\begin{corollary}\label{nearisom} $G$ is near-isometric to $S/P$ with constant $c=1$.

\end{corollary}

\subsection{Definition of the language}
In this section we define the regular language we will use to show $G$ has rational growth.  In order to apply Theorem \ref{PutFFT}, we need to construct a regular language and an acceptor for the partition $P$ defined above. Having an acceptor for the partition means being able to tell when two words belong to the same equivalence class modulo $P$.  In light of Proposition \ref{grpequal}, this implies that we must be able to simulate Laurent polynomial long division in an FSA. Because an FSA can only store a finite amount of memory at each stage of computation, we need a finite criterion which distinguishes words in different equivalence classes.  To this end we will show that minimal length words have bounded coefficients, and that after polynomial long division, coefficients stay bounded.  

\subsubsection{Bounding Coefficients}
The following lemma bounds coefficients of minimal words.  Here it will be clear why the conditions guaranteed by Lemma 1 are important.  Recall that $char(A)=p_1\cdots p_k$ where $\deg(p_i)=d_i$, and $M_i$ denotes the maximal coefficient of $p_i$.

\begin{lemma}\label{bound1}
Let $(\overl{t},h)\in S$ be a word with minimal length in its equivalence class modulo $P$, where $\lp{t}{k}$. Then for each $i=1,\ldots, k$, $|t_{i}|_{\infty}<2N|M_i|$.
\end{lemma}
\footnotesize{PROOF}: 
\normalsize We prove the contrapositive.  Suppose an unreduced type $\overl{t}=(t_1(z_1),\ldots, t_k(z_k))$ has some coefficient $c_{ij}$ with $|c_{ij}|>2N|M_i|$.  Note that by adding or subtracting a term of the form $(0,\ldots, z_{i}^rp_i,\ldots,0)$, we do not change the type of $\overl{t}$, since the difference is divisible by $(p_1,\ldots, p_k)$. Choose $r$ so that the maximal coefficient of $z_{i}^rp_i$ multiplies $z_i^j$.  Then we can reduce $|c_{ij}|$ by $2NM_i$ while increasing the other coefficients by a maximum of 
\begin{align*} 2(d_i-1) +2N(M_i-1)&\leq 2(N-1) +2N(M_i-1)\\ &<2NM_i.\end{align*} Thus $\overl{t}$ was not minimal.  \QED

Note that the previous lemma implies that, in particular, $|\overl{t}|_{\infty} < 2N\max\{|M_i|\}$ for minimal unreduced types $\overl{t}$.  The next lemma states that if $(p_1,\ldots,p_k)$ divides the difference of two unreduced types, the coefficients of the quotient remain bounded. 

\begin{lemma}\label{bound2}  

For every integer $A>0$, there exists $B_A>0$ so that if \[|\overl{t}_1|_{\infty},|\overl{t}_2|_{\infty} \leq A\] and if $(p_1,\ldots,p_k)$ divides $\overl{t}_1-\overl{t}_2$, then \[|(\overl{t}_1-\overl{t}_2)/(p_1,\ldots,p_k)|_{\infty}\leq B_A.\]

\end{lemma}
\footnotesize{PROOF}: 
\normalsize Observe that because $2|M_i|>|p_i|_{1}$, the largest coefficient of $g(z_i)p_i(z_i)$ is at least as large as that of $g$. If $\overl{t}_1-\overl{t}_2 = (g_1p_1,\ldots,g_kp_k)$, then by assumption $|\overl{t}_1-\overl{t}_2|_{\infty}\leq 2A$, hence $|(g_1,\ldots,g_k)|_{\infty}\leq2A$. Setting $B=2A$ suffices.  \QED

This final lemma guarantees that during the process of polynomial long division, coefficients of the remainder do not get too large either.    

\begin{lemma}\label{divide} Let $A, B >0$.  $\exists C(A,B)>0$ so that the following holds:
		
Given  $\overl{t}_1$, $\overl{t}_2$, suppose that \[|\overl{t}_1|_{\infty},|\overl{t}_2|_{\infty} \leq A\] and that for each $i=1,\ldots,k$ we can write \[(\overl{t}_1-\overl{t}_2)_i = p_i(z_i)q_i(z_i) + r_i(z_i)\] where $r_i(z_i)$ is a polynomial of degree at most $d_i-1$ and  $|q_i|_{\infty}<B$.  Then $|r_i|_{\infty}<C(A,B)$.
		
\end{lemma}
\footnotesize{PROOF}: 
\normalsize The coefficients of $p_i(z_i)q_i(z_i)$ are bounded by $B(|p_i|_{1})<2BM_i$, while the coefficients of  $\overl{t}_1-\overl{t}_2$ are bounded by  $2A$.   Then $|r_i|_{\infty}< 2A + 2BM_i$.  Choosing $C(A,B) = 2BM$, where $M=\max\{M_i\}$ proves the lemma.\QED 

Now we are ready to define the language.  Let $A_n= \{-n,\ldots, n\}^k \times \{-1,1,2\}$ be an alphabet with weighting \[\phi((a_1,\ldots,a_k),b)=\sum_{i=1}^k |a_i|+|b|.\]  

We define languages $\Lambda_n, \Lambda_n'$ following Putman. Words in $\Lambda_n$ are exactly those of the form: \[(\cdot,2)\cdots(\cdot,2)(\cdot,\pm1)\cdots(\cdot,\pm1)(\cdot,2)\cdots(\cdot,2)\]  satisfying 
\begin{enumerate}[(1)]
	\item There is at least one center $(\cdot,\pm1)$ term.
	\item All center terms all have the same sign in the second coordinate. 
	\item If the second coordinate of the center terms is $-1$, there are at least two center terms.
	\item If the first or last letters are of the form $((a_1,\ldots,a_k),2)$ then the $a_i$ are not all 0. 
	
\end{enumerate}
The tripartite division of each word reflects the tail, center, head division of unreduced types, while the sign of the second coordinate in the middle terms tells the sign of the height.  $\Lambda_n'$ is defined in the same way but without condition (4).

It is clear that each of these languages is regular, for each can be written as a regular expression.  Finally we define a map $\psi':\Lambda_n'\rightarrow \mlaur$ by 
\begin{align*} \psi'&\left(\prod_{i=1}^{n_1}((\right.c_{1i},
\ldots,c_{ki}),2)\prod_{i=0}^{n_2}((c_{1i}',
\ldots,c_{ki}'),+1)\left.\prod_{i=1}^{n_3}((c_{1i}'',
\ldots,c_{ki}''),2)\right) =\\
&\sum_{i=1}^{n_1}(c_{1i}z_1^{i-n_1-1},\ldots,c_{ki}z_k^{i-n_1-1})+\sum_{i=0}^{n_2}(c_{1i}'z_{1}^{i},\ldots,c_{ki}'z_{k}^{i})+\sum_{i=1}^{n_3}(c_{1i}''z_1^{i+n_3},\ldots,c_{ki}''z_k^{i+n_3}),
\end{align*}

and

\begin{align*} \psi'& \left( \prod_{i=1}^{n_1}((c_{1i},
 \ldots,c_{ki}),2) \prod_{i=0}^{n_2}((c_{1i}',
\ldots,c_{ki}'),-1) \prod_{i=1}^{n_3}((c_{1i}'',
\ldots,c_{ki}''),2)\right) =\\ 
&\sum_{i=1}^{n_1}(c_{1i}z_1^{i-n_1-n_2-1},\ldots,c_{ki}z_1^{i-n_1-n_2-1})+\sum_{i=0}^{n_2}(c_{1i}'z_1^{i-n_2},\ldots,c_{ki}'z_k^{i-n_2})+\sum_{i=1}^{n_3}(c_{1i}''z_1^{i},\ldots,c_{ki}''z_k^{i}).
\end{align*}

We extend this to $\Psi': \Lambda_n\rightarrow S$ in the following way.  Let $w \in \Lambda_n$ and let $h=$ number of central $(\cdot, \pm1)$ terms in $w$.  Then

\[\Psi'(w) = (\psi'(w),\pm h).\] 
We obtain a map $\Psi: \Lambda_n\subset\Lambda_n'\rightarrow S$  by restriction.  Conditions (3) and (4) imply that $\Psi$ is injective.  The partition $P$ on $S$ induces partitions $\Pi_n',\Pi_n$ on $\Lambda_n', \Lambda_n$ respectively, and induced maps $\overl{\Psi}':\Lambda_n'/\Pi_n'\rightarrow S/P$,  $\overl{\Psi}:\Lambda_n/\Pi_n\rightarrow S/P$. Lemma \ref{bound1} implies that when $n\geq2N\max\{|M_i|\}$, $\overl{\Psi}'$ and $\overl{\Psi}$ are surjective, hence $\overl{\Psi}$ is bijective. This discussion now implies:

\begin{thm}\label{isom} 
For $n\geq M=2N\max\{|M_i|\}$, the induced map $\overl{\Psi}:\Lambda_n/\Pi_n\rightarrow S/P$ is an isometry.  
\end{thm} 

\subsection{The Acceptor Automaton}
Using the map $\Psi': \Lambda_n'\rightarrow S$ we can associate to a word $w \in \Lambda_n'$ the tuples of Laurent polynomials \[Tail(w) = \Tl{\Psi'(w)},\]  \[Center(w) = \Ct{\Psi'(w)},\] \[Head(w) = \Hd{\Psi'(w)}\] We also define analogues of the length functions $T(w)$ (resp.  $H(w))$ to be the maximal number of $(\cdot,2)$ letters in any $w_i$ before (resp. after) the middle $(\cdot,\pm1)$ letters.  Thus even letters of the form $((0,\ldots,0),2)$ contribute to the length of the tail or head.  Note that if $w\in\Lambda_n'\setminus \Lambda_n$ then either $T(w)\neq T_h(\Psi'(w))$ or $H(w)\neq H_h(\Psi'(w))$.  

As above with $\Lambda_n$, we will define a sequence of padded languages indexed by $n$: \[R_{n,i}=\left\{\begin{array}{c|c} &\overl{w}_1=\overl{w}_2 \text{ modulo }\Pi_n,\\ (w_1,w_2)\in \Lambda_n\times\Lambda_n & |T_h(w_1)-T_h(w_2)|\leq i,\text{ and }\\&|H_h(w_1)-H_h(w_2)|\leq i\end{array}\right\}.\]

The main result of this section is the following:

\begin{thm}\label{acceptor} $R_{n,i}$ is an acceptor for the partition $\Pi_n$ of $\Lambda_n$.   
\end{thm}
We remark that the main work in proving the theorem will be to show that $R_{n,i}$ is regular.  It follows by definition of the language $R_{n,i}$ that it is an acceptor for the partition $\Pi_n$.  In order to prove regularity we will need the next lemma, which defines an auxiliary language $R_n'$:

\begin{lemma}\label{regularacc} Define the padded language
\[R_{n}'=\{(w_1,w_2)\in \Lambda_n'\times\Lambda_n'| \overl{w}_1=\overl{w}_2 \text{ modulo }\Pi_n' \text{, } T(w_1)=T(w_2),  \text{ and } H(w_1)=H(w_2)\}.\] Then $R_n'$ is  regular.  
\end{lemma}
\footnotesize{PROOF}: 
\normalsize 

We will show that rev$(R_n')$ is regular, by Theorem \ref{ECHprop} this is sufficient.  The automaton we construct will read $w_1$ and $w_2$ from right to left, and simulate the division of $w_1-w_2$ by $(p_1,\ldots,p_k)$ at each step, keeping track of the remainder.  By Lemmas \ref{bound2} and \ref{divide}, we know that if the coefficients of the remainder ever exceed $C(n, B_n)$, then $w_1\neq w_2$ modulo $\Pi_n$.  The automaton also keeps track of the head, center and tail of $w_1$ and $w_2$ to ensure that they line up.  

Let $B= B_n$ and $C = C(n,B_n)$ as in Lemmas  \ref{bound2} and \ref{divide}.  The automaton has a fail state and the following additional states:\[\left\{\begin{array}{c|c}

&r_i = c_{i0}+c_{i1}z_i+\cdots+c_{id_i}z_i^{d_i-1}\mbox{ with } c_{ij}\in \Z, |c_{ij}|\leq C,\\
((r_1,\ldots,r_k),p)& \mbox{ and }\\& p \in \{\mathcal{H},\mathcal{T},\mathcal{C}_1,\mathcal{C}_{-1},\mathcal{C}_{-1,1}\}
\end{array}\right\}.\]
The $r_i$ represent the remainder in the $i$th component, and $p$ describes what part of word we are currently reading.  The start state is $((0,\ldots,0),\mathcal{H})$ since we are reading from right to left.  Let $w_1=\Psi'^{-1}((t_{11},\ldots,t_{1k}), h_1)$ and $w_2=\Psi'^{-1}((t_{21},\ldots,t_{2k}), h_2)$ We describe how the computation proceeds in each coordinate.  If in the $j$th coordinate, after reading $l$ letters we are in state $(r_j,p)$ then we have:\[(t_{1j}-t_{2j})_{N_2-l+1}^{N_2}= z_j^{N_2-l+1}(q(z_j)p_j(z_j)+r_j).\]
Here $q(z_j)$ is a Laurent polynomial with coefficients bounded by $B$, whose actual value does not matter. We examine the effect of adding the next term, which we represent in the form $z_j^{N_2-l}(a_{1}-a_{2})$:
\begin{align*}(t_{1j}-t_{2j})_{N_2-l}^{N_2}&= z_j^{N_2-l+1}(q(z_j)p_j(z_j)+r_j)+z_j^{N_2-l}(a_{1}-a_{2})\\
&=z_j^{N_2-l}(z_jq(z_j)p_j(z_j)+z_jr_j+(a_{1}-a_{2})).
\end{align*}
Note that $z_jr_j+(a_{1}-a_{2})$ has degree at most $d_j$. If $z_jr_j+(a_{1}-a_{2})$ has degree $d_j$ we divide it by $p_j$ to obtain $r_j'$, otherwise $r_j'=z_jr_j+(a_{1}-a_{2}))$. If at any stage, this division procedure yields an $r'_j$ whose coefficients have absolute value greater than $C$, we fail.  In this way we can compute the remainder in each component based on the input and the previous remainder.  Because the coefficients are bounded, these transitions are all determined in advance.  

In the second coordinate, we begin in state $\mathcal{H}$ and stay in that state as long as we read $(\cdot,2)$.  If we are in state $\mathcal{H}$ and read $(\cdot,1)$ in both $w_1$ and $w_2$ at the same time, we transition to state $\mathcal{C}_1$. If we are in state $\mathcal{H}$ and read $(\cdot,-1)$ in both $w_1$ and $w_2$ at the same time, we transition to state $\mathcal{C}_{-1,1}$.  If, while in state $\mathcal{C}_{-1,1}$, if we read $(\cdot,-1)$ in both $w_1$ and $w_2$ again, we transition to state $\mathcal{C}_{-1}$, otherwise we fail.  This is to ensure that both $w_1$ and $w_2$ have the form required by condition (3) in the definition of $\Lambda_n'$.  Finally, if we are in either of the states $\mathcal{C}_{\pm1}$ and we read $(\cdot,2)$ in both $w_1$ and $w_2$ at the same time, we transition to state $\mathcal{T}$.  If at any point, the second entries in both $w_1$ and $w_2$ don't match up, or proceed in an order different from $2\rightarrow\pm1\rightarrow2$, we fail.  The accept states are $((0,\ldots,0),\mathcal{T})$, $((0,\ldots,0),\mathcal{C}_1)$ and $((0,\ldots,0),\mathcal{C}_{-1})$.  It is clear that we end in one of these states iff $w_1$ and $w_2$ are in the form required by $\Lambda_n'$ and their difference is divisible by $(p_1,\ldots, p_k)$.\QED

We now deduce the theorem from the lemma:

\footnotesize{PROOF} 
\normalsize of Theorem \ref{acceptor}:

We represent $R_{n,i}$ as a union of languages, each of which is regular by the previous lemma and Theorem \ref{ECHprop}.  Let $1\leq j\leq i$.  Set \[Q_{j}= \prod_{k=1}^{j}((0,\ldots,0),2).\]
Then \[
R_{n,i} = \bigcup_{1\leq j_1,j_2\leq i} \left\{\begin{array}{c|c} &  \text{Either }
	(Q_{j_1}w_1Q_{j_2},w_2)\text{ or }\\ (w_1,w_2)\in  \Lambda_n\times \Lambda_n& (w_1Q_{j_1},w_2Q_{j_2})\text{ or }
	(w_1Q_{j_1},Q_{j_2}w_2)\text{ or }\\&(w_1,Q_{j_1}w_2Q_{j_2}) \in R_n' \end{array} \right\}.\]
Each of the languages in the union is constructed from a regular language by first-order predicates and concatenation, and hence is regular by Theorem 3.  Since there are only a finite number of affixes $Q_j$, we have represented $R_{n,i}$ as a finite union of regular languages.\QED

\section{FFTP and the proof of the main theorem}\label{FFTP}

This section is devoted to the proof of the FFTP for the pairs $(\Lambda_n, \Pi_n)$ for $n$ sufficiently large.  Together with Theorem 5 this implies rational growth of the group $G$.  First we need a definition which serves as  a kind of Hausdorff distance for types and heights: 

\begin{definition} Let $f_1=\sum_ic_{1i}z^i$ and $f_2=\sum_ic_{2i}z^i$ be two Laurent polynomials in $z$.    The \textbf{divergence} of $f_1$ and $f_2$ is defined as \[\dv(f_1,f_2)=
\max_n\sum_{i=-\infty}^n|c_{1i}|-|c_{2i}|\]
as $n$ ranges from $-\infty$ to $\infty$.  If $\overl{t}=(t_1,\ldots, t_k)$ and $\overl{s}=(s_1,\ldots, s_k)$ are two tuples of Laurent polynomials then $\Dv(\overl{t},\overl{s})$ is defined to be $\displaystyle \sum_{j=1}^k |\dv(t_j,s_j)|$, and if $w_1, w_2\in \Lambda_n$ are words then $\Dv(w_1,w_2) = \Dv(\utp{(w_1)},\utp{(w_2}))$.
\end{definition}

The theorem will follow from the following key lemma.  

	\begin{lemma}\label{fellowtravel} Set $M=\max_i\{|M_i|\}$.  $\exists$ constants $K_1,K_2,K_3$ such that $\displaystyle K_1\geq 2NM$  and
	\begin{enumerate}
		\item If $w_1\in \Lambda_{2NM}$ is not minimal modulo $\Pi_{2MN}$, $\exists w_2\in \Lambda_{K_1}$ so that 
			\begin{itemize}
				\item $w_1=_Gw_2$ and $\norm{w_1}>\norm{w_2}$.
				\item $|T(w_1)-T(w_2)|,|H(w_1)-H(w_2)|\leq K_2$.
				\item $\Dv(w_1,w_2)\leq K_3$. 
			\end{itemize}
		\item If $w_1=_Gw_2$ are two different minimal length representatives modulo $\Pi_M$, then $|T(w_1)-T(w_2)|,|H(w_1)-H(w_2)|\leq K_2$.
	\end{enumerate}
	\end{lemma}
	
Before we go into the proof of the lemma, we deduce the main theorem:

\footnotesize{PROOF} 
\normalsize of Theorem \ref{main1}:  

By Corollary \ref{nearisom} and Theorems \ref{PutFFT} and \ref{isom}, it suffices to show that for $n$ sufficiently large $\Lambda_n/\Pi_n$ has the FFTP. Let $K_1$, $K_2$ and $K_3$ be the constants from Lemma \ref{fellowtravel}.  We demonstrate the FFTP for the pair $(\Lambda_{2MN},\Pi_{2MN})$ with acceptor $R_{K_1,K_2}$ with  fellow-traveling constant $K_3 +(kK_1+6)K_2$.

First note that since $K_1\geq 2NM$, $\overl{\Psi}$ is surjective.  Now consider some $w_1\in \Lambda_{2MN}$ which is not minimal. Lemma 13 guarantees there exists $w_2\in \Lambda_{K_1}$ satisfying the first two conditions in the definition of the FFTP.  Hence all we need to show is that $w_1$ and $w_2$ are fellow-travelers with constant  $K_3 +(k K_1+6)K_2$. 

Without loss of generality, we assume that $ T(w_1)\geq T(w_2)$; the proof in the reverse case is similar. The difference in length between the tail of $w_1$ and that of $w_2$ is at most $K_2$, and $\dv(w_1,w_2)\leq K_3$.  Let $v_1$ and $v_2$ denote length $i$-initial segments of $w_1$ and $w_2$ respectively, and consider the difference \[|\norm{v_1}-\norm{v_2}|.\]  The portion of $v_1$ beyond the tail of $v_2$ contributes at most $2K_2$ to the difference, and similarly a portion of the head of $v_1$ or $v_2$ beyond that of the other contributes at most $2K_2$.  Where they overlap, the contribution of the divergence is at most $K_3$.  Finally, the remaining portion of $v_2$ where the two words do not overlap has length at most $K_2$ and each term contributes at most $k K_1 +2$ to the length difference.  Thus we obtain:

\[|\norm{v_1}-\norm{v_2}|\leq 4K_2+K_3 +(k K_1+2)K_2=K_3 +(k K_1+6)K_2.\]

The second part of Lemma \ref{fellowtravel} ensures that if $w_1$ and $w_2$ are both minimal, then $(w_1, w_2) \in R_{K_1,K_2}$.  This proves all the conditions in the FFTP, hence the theorem.  \QED

We now proceed with the proof of the lemma.

\footnotesize{PROOF} 
\normalsize of Lemma \ref{fellowtravel}:  

Let $B = B_{2MN}$ be the constant from Lemma \ref{bound2}.  We will show that the following constants satisfy the claims of the lemma: 
\begin{itemize}
\item $K_1 = 2MN(B+2)$,
\item $K_2 = 2MNB$,
\item $K_3 = k\left(k(12MNB+4MN+1)+2MB + 2MNB+4M\right).$
\end{itemize}

Given a word $w_1\in \Lambda_{2MN}$ which is not minimal modulo $\Pi_{2MN}$, there exists some representative $w_2\in  \Lambda_{2MN}$, with $\overl{w}_1=_G\overl{w}_2$ and $\norm{w_1}>\norm{w_2}$.  Writing $\Psi(w_i)=(\overl{t}_i, h)$ we know that by Lemma \ref{bound2} there exists a tuple of Laurent polynomials $\overl{q}=(q_1,\ldots, q_k)$ satisfying $\overl{t}_2 = \overl{t}_1+ (p_1q_1,\ldots,p_kq_k)$ and $|\overl{q}|_{\infty}<B$.  For each $i$, $1\leq i \leq k$, let $t_{1i}$ denote the $i$th component of $\overl{t}_1$, $t_{2i}$ denote the $i$th component of $\overl{t}_2$ and $q_{i}$ denote the $i$th component of $\overl{q}$ so that we have
\[t_{2i}(z_i)= t_{1i}(z_i)+p_i(z_i)q_i(z_i).\]
Recall that the degree of $p_i$ is $d_i$, and that the magnitude of largest coefficient of $p_i$ is $|M_i|$.

The strategy of the proof is as follows.  We will modify $\overl{t}_2$ and $\overl{q}$ to produce a new word $w_2'$ which is \emph{closer} to $w_1$ but shorter.  However, in order to produce such a word, we may need to allow coefficients in $\utp(w_2')$ to be as large as $K_1$.  More precisely, we will modify $\overl{t}_2\Rightarrow \overl{t}_2'$ and $\overl{q}\Rightarrow \overl{q}'$ so that by setting $w_2'= \Psi^{-1}(\overl{t}_2',h)$ we have:
\begin{enumerate}[(1)]
	\item $|\overl{q}'|_{\infty}\leq B$ and $\norm{w_2'}<\norm{w_2}$.  
	\item $T(w_2')-T(w_1)\leq K_2$ and  $H(w_2')-H(w_1)\leq K_2$.
	\item $T(w_1)-T(w_2')\leq K_2$ and $H(w_1)-H(w_2')\leq K_2$.
	\item For each $i$, $1\leq i \leq k$, $\dv(t_{2i}',t_{1i}) \leq K_3/k$.
	\item For each $i$, $1\leq i \leq k$, $\dv(t_{1i},t_{2i}') \leq K_3/k$.
\end{enumerate}
Note that the first part of (1) implies that $w_2'\in \Lambda_{K_1}$.

The proof will proceed in 4 steps, each step guaranteeing that some collection of conditions (1)--(5) above can be achieved for the tuple of Laurent polynomials.  Each step will provide a construction which may be applied to a single Laurent polynomial.  After defining how the construction is carried out on each component separately, we will indicate how to extend the construction to the tuple.  Because each construction only applies to a single Laurent polynomial, we will greatly simplify our exposition by allowing the following notation to respresent any given component of the Laurent polynomial:

\begin{itemize}
	\item $\displaystyle t_{1}(z)=\sum_{i=-\infty}^{\infty}a_{i}z^i$, $\displaystyle t_{2}(z)=\sum_{i=-\infty}^{\infty}b_{i}z^i$ and  $\displaystyle q(z) = \sum_{i=-\infty}^{\infty}c_{i}z^i$.
	\item $p(z)$ (corresponding to $p_i(z_i))$ will denote an irreducible characteristic polynomial of degree $d$ $(=d_i)$.
	\item $N$, $M$, and $B$ will all remain as above, even though they figure into the proof of the special case of one polynomial below.  
 \end{itemize}

For ease of notation, in each stage of the modification we will reset the notation above -- $w_1$ will refer to the original word, $w_2$ the minimal length word, and $w_2'$ the modification.  

\textbf{Step 1:} We can alter $w_2$ so that (1) and (2) are satisfied.  

We will only demonstrate how to achieve $T(w_2')-T(w_1)\leq K_2$, the corresponding inequality for the heads is similar.  Assume that $T(t_2)>T(t_1)+ 2MNB$.  In this case the lowest degree of any non-zero term of $q$ must be smaller than that of $t_1$.  We will find a word with shorter tail satisfying condition (1).  Iterating this procedure will complete this part of the modification.  

Let $D$ be the smallest index for which $q$ has a non-zero coefficient. Thus  we have
\[t_{2} = t_{1} +p(z)\sum_{i=D}^{\infty}c_{i}z^i.\]

Set $D'= D+2MNB$, and define 
\[t_{2}' = t_{1}+p(z)\sum_{i=D'}^{\infty}c_{i}z^{i}.\] 
By construction we have reduced the difference between the lengths of tail of $t_1$ and $t_2$ by at least $2NMB$.  The discrepancy may be strictly greater in the case that $c_{D'}=0$, for example. Also, note that we have the following inequalities \[\left\{\begin{array}{ll}|b_{i}'|  =0 & \mbox{ if } i\leq D'\\
		 |b_{i}'| \leq |b_{i}|+2MB& \mbox{if }  D'\leq i\leq D'+d\\
		 |b_{i}'| =  |b_{i}|& \mbox{ if } j\geq D'+d
\end{array}\right.\]

We may have  increased the length of the head $t_2$ by at most $d-1$:\[H(t_2)-H(t_2')\geq 1-d.\]  Apply the above construction to each component $i$ which satisfies $T(t_{2i})>T(t_{1i})+ 2MNB$, and set $w_2'=\Psi^{-1}(t_{2}',h)$. Combining these estimates  for each $i$ we have: \begin{align*}
\norm{w_2}-\norm{w_2'}&=2(T(w_2)-T(w_2')) +2(H(w_2)-H(w_2')) +\sum_{i=1}^{k}|t_{2i}|_{1}-|t_{2i}'|_{1}\\
				&\geq 2(2MNB)-2(\max_i\{d_i\}-1)-2MB\sum_{i=1}^kd_i\\
				&\geq4MNB-2(N-1)-2MNB\\
				&> 2N(MB-1)>0
\end{align*} 
since $B>M>1$ and $N=\sum_i d_i$. \QED

\textbf{Step 2:} We can alter $w_2$  from the previous step so that (1)--(3) are satisfied.

Assume that $T(t_{1})>T(t_{2})+2MNB$. In this case, the lowest degree non-zero coefficient of $q$ is at least as small as that of $t_1$, but when added together, some cancellation must occur. We can modify $w_2$ by deleting the top of $q$,  and thereby ensuring fewer terms in the tail of  $t_2$ cancel out.   

Now let $D$ denote the index of the lowest degree non-zero term of $t_{1}$.  Set \[D'=D+2MNB\]
 and define
\[t_{2}' =t_{1}+p(z)\sum_{i=-\infty}^{D'-1}c_{i}z^{i}.\] 


We want to ensure that the coefficient $b_{D}'$ is non-zero.  Thus we may need to add $\pm z^{D}p(z)$ in such a way that $|t_2'|_{\infty}\leq K_1$.  By construction we have
\[T(t_1)-T(t_2') = 2MNB. \]
But the possible addition of $\pm z^{D}p(z)$ means the head of $t_2'$ may differ from that of $t_1$ by at most $d-1$:
\[|H(t_1)-H(t_2')|\leq d-1< 2MNB.\]
These two observations imply that the requirements of conditions (2) and (3) have been met.  
Tracing through the modifications, we obtain the following inequalities:

\[\left\{\begin{array}{ll}|b_{i}'|  =0 & \mbox{if } i< D'\\
		 |b_{i}'| \leq |a_{i}|+2M(B+1)& \mbox{if }  D'\leq i\leq D'+d-1\\
		 |b_{i}'| =  |a_{i}|& \mbox{if } i\geq D'+d 
\end{array}\right.\]

As in the previous step, we apply the above construction to each component of the tuple of Laurent polynomials where $T(t_{1i})>T(t_{2i})+ 2MNB$.  Set $w_2'=\Psi^{-1}(t_{2}',h)$.  Because conditions (2) and (3) are satisfied component-wise, they are satisfied for $w_2'$.  Now we need only make sure that $\norm{w_2'}<\norm{w_1}$ so that condition (1) is still satisfied. Combining these inequalities we obtain:
\begin{align*}
\norm{w_1}-\norm{w_2'}&=2(T(w_1)-T(w_2')) +2(H(w_1)-H(w_2')) +\sum_{i=1}^{k}|t_{1i}|_1-|t_{2i}'|\\
				&\geq 2(2MNB)-2(\max_i\{d_i\}-1)-2M(B+1)\sum_{i=1}^kd_i\\
				&\geq4MNB-2(N-1)-2M(B+1)N=2MN(B-1)-2(N-1)\\
				&> 2N(M(B-1)-1)>0
\end{align*} 
since $B>M>1$ and $N=\sum_i d_i$. Note that in the above calculation we have ignored the contribution of all coefficients when $i<D'$.   The modification in the corresponding case for heads is similar.  \QED

\textbf{Step 3:} We can modify $w_2$ from Step 2 so that conditions (1)--(4) are all satisfied.  

Suppose that for some component and some minimal $D\in \Z$ we have \[\sum_{i=-\infty}^D|b_i|-|a_i|> 10MNB+4M+1\] so that in particular $\Dv(w_1,w_2)\geq10MNB+4M+ 1$.  We will now produce a $t_2'$ satisfying conditions (1)--(3).  By iterating this procedure we will be able to reduce the divergence to be at most $10MNB+4M+1$. 

The idea is simply to delete the first portion of $w_2$ which makes the divergence larger than $10MNB+4M+1$.  First define an auxiliary Laurent polynomial $t_2''$ as 
\[t_2'' = t_1+\sum_{i=D+1}^{\infty}c_iz^i\]
We obtain the following estimates:
\[\left\{\begin{array}{ll}|b_{i}''|  =a_i & \mbox{if } i< D+1\\
		 |b_{i}''| \leq |b_{i}|+2MB & \mbox{if }  D+1\leq i\leq D+d\\
		 |b_{i}''| =  |b_{i}|& \mbox{if } i\geq D+d+1
\end{array}\right.\]
As in the previous step, we may have deleted too much from either the head or the tail.  We therefore form a new Laurent polynomial $t_2'$ by adding $p(z)(A_1z^{D_1} + A_2z^{D_2})$ where the $A_i\in\{0,\pm1\}$, $D_i\in \Z$ and are chosen in such a way that \[|T_h(t_2')-T_h(t_1)|\leq2MNB \mbox{ and } |H_h(t_2')-H_h(t_1)|\leq2MNB.\] Lastly, observe that the divergence of $t_2''$ and $t_2'$ is at most $4M$.

We now apply this construction to each component where \[\sum_{i=-\infty}^D|b_i|-|a_i|> 10MNB+4M+1.\]  Then if we set $w_2'=\Psi^{-1}(t_2',h)$ we need only check that $\norm{w_2'}<\norm{w_2}$, since by repeating this process we can obtain $\dv(t_{1i},t_{2i}')\leq 10MNB+4M+ 1$. Suppose that we need to apply the construction to $k'$ indices, where $1\leq k'\leq k$. From the construction  and the above estimates we have that:
\begin{align*}
\norm{w_2}-\norm{w_2'}&=2(T(w_2)-T(w_2')) +2(H(w_2)-H(w_2')) +\sum_{j=1}^{k'}|t_{2i_j}|_1-|t_{2i_j}'|_1\\
				&\geq 2(-2MNB)+2(-2MNB) +\sum_{j=1}^{k'}|b_{i_j}|-|b_{i_j}''|+\sum_{i=1}^{k'}|b_{i_j}''|-|b_{i_j}'|\\
				&\geq-8MNB +k'(10MNB+4M+ 1-2MNB)-4Mk' \\
				&\geq (k'-1)(8MNB)+k'>0.
\end{align*} 
Since there is at least one component which we need to modify, $k'\geq 1$.  To obtain the second line, we used the fact that $|T(w_2)-T(w_1)|, |H(w_2)-H(w_1)|\leq 2MNB$. Note that we do not consider the contribution from those indices we did not alter, since on these, $w_2$ and $w_2'$ agree.   \QED

\textbf{Step 4:} We can modify $w_2$ from Step 3 so that conditions (1)--(5) are all satisfied.

 For this case we will present a general modification which will depend on how large the divergence is in each component.   Assume that for some component there exists a $D\in \Z$ and a positive integer $n_0$ so that \[ \sum_{i=-\infty}^D|a_i|-|b_i|> n_0(12MNB+4MN+1).\]  Choose $D_0$ minimal with this property.  Since for any $i$, $||a_i|-|b_i||\leq 2MB$ we have
\[\sum_{i=-\infty}^{D_0}|a_i|-|b_i|\leq n_0(12MNB+4MN+1)+2MB\] by the minimality of $D_0$. In this case we delete terms in $q$ past $D_0$ and hope to show we have reduced the divergence.  Define $t_2''$ by \[t_2'' = t_1+p(z)\sum_{i=\infty}^{D_0}c_iz^i.\]  We can bound the coefficients of $t_2''$ as follows:
\[\left\{\begin{array}{ll}|b_{i}''|  =|b_i| & \mbox{if } i\leq D_0\\
		 ||b_{i}''| - |b_{i}||\leq 2MB & \mbox{if }  D_0+1\leq i\leq D_0+d\\
		 |b_{i}''| =  |a_{i}|& \mbox{if } i\geq D_0+d+1
\end{array}\right.\]
As in Step 3, we compensate for the fact that we may have lopped off too much of $t_2$ by adding $p(z)(A_1z^{D_1} + A_2z^{D_2})$  to $t_2''$ where the $A_i\in\{0,\pm1\}$, $D_i\in \Z$. This adjusted Laurent polynomial we name $t_2'$.  The divergence of $t_2''$ and $t_2'$ is at most $4M$ and by the minimality of $D_0$ the divergence of $t_2'$ and $t_1$ is less than \[n_0(12MNB+4MN+1)+2MB + 2MNB+4M.\]
In order to prove the result, we will consider several cases.  The problem is that although the contribution of one component to the divergence may be very large, if we delete it, we need to make sure that we have not made the length of $w_2'$ greater than $w_1$.  Thus we distinguish several cases where some collection of the components have large divergence, and the other components have divergence bounded by $10MNB+4M+1$. Thus we define cases $\textbf{C}_n$, where $0\leq n\leq k-1$:
\[
 \textbf{C}_n:\begin{array}{c}  \mbox{ For $n$ components, $\dv(t_{1i},t_{2i})\leq12MNB+4MN+1$ and}\\ \mbox{ for $k-n$ components, } \dv(t_{1i},t_{2i})>(n+1)(12MNB+4MN+1).
\end{array}\]
Although these cases are not exhaustive of every possibility, their complement in $(\R_{\geq0})^k$ is contained in the bounded region defined by the $k$ inequalities \[|\dv(t_{1i},t_{2i})|\leq k(12MNB+4MN+1)+2MB + 2MNB+4M.\]

If the pair $\overl{t}_1,\overl{t}_2$ falls into the case $\textbf{C}_n$, apply the above construction with constant $n_0=n+1$ to each of the $k-n$ components where necessary and set $w_2'=\Psi^{-1}(t_2',h)$. The final step is to show that $\norm{w_2'}<\norm{w_1}$.  Without loss of generality, we assume that we have modified the last $k-n$ components. The other $n$ have not changed, hence from the calculations above:
\begin{align*}
\norm{w_1}-\norm{w_2'}&=2(T(w_1)-T(w_2')) +2(H(w_1)-H(w_2')) +\sum_{i=1}^{n}|t_{1i}|_1-|t_{2i}'|_1+\sum_{i=n+1}^{k}|t_{1i}|_1-|t_{2i}'|_1\\
				&\geq 2(-2MNB)+2(-2MNB) -\sum_{i=1}^{n}(12MNB+4MN+1)\\&\hspace{50pt}+\sum_{i=n+1}^{k}|t_{1i}|_1-|t_{2i}''|_1+\sum_{i=n+1}^{k}|t_{2i}''|_1-|t_{2i}'|_1\\
				&\geq-8MNB -n(12MNB+4MN+1)+(k-n)(n+1)(12MNB+4MN+1)\\
				&\hspace{50pt}-(k-n)2MB-(k-n)(4M+2MB).\\
\end{align*} 
Since $1\leq k-n<N$, the last inequality is greater than \[-12MNB+(12MNB+4MN+1)-4MN\geq1>0.\]
By the above discussion, we can set \[K_3=k(k(12MNB+4MN+1)+2MB + 2MNB+4M)\] since the divergence of $w_1$ and $w_2$ is the sum of the divergences of each component.
Note that $K_3>10MNB+4M+1$ so the same constant works for both Steps 3 and 4.  This completes the proof of Step 4.\QED

Going through each step of the proof, we see that since we had to add at most 2 constant multiples of polynomials at any stage, we can choose \[K_1 = 2MN(B+2) .\]

Steps 1 and 2 imply that we can choose
\[K_2 = 2MNB,\]
and the proof of Step 4 implies that as an upper bound on the divergence we may take
\[K_3 = k(k(12MNB+4MN+1)+2MB + 2MNB+4M).\] This concludes the proof of the first part of Lemma \ref{fellowtravel}.

For the second part, if $w_1$ and $w_2$ are equal modulo $\Pi_{2MN}$, but we have either $|T(w_1)-T(w_2)|>K_2$ or $|H(w_1)-H(w_2)|>K_2$, then Steps 1 and 2 imply we can find $w_2'\in \Lambda_{K_1} $ with $w_2'=w_1$ modulo $\Pi_{K_1}$ and $\norm{w_2'}<\norm{w_1}$.  This is impossible, since $w_1$ was assumed minimal.  \QED

\begin{center}\Large \textbf{Part II: Almost Convexity}
\end{center}

\section{Almost Convexity}\label{AC}
Let $G$ be a finitely presented group with generating set $S=S^{-1}$.  Equip $G$ with the word metric $d_S$ coming from $S$ and denote by $\beta(n)$ and $\Sigma(n)$ the ball and sphere of radius $n$ respectively, centered on the identity in the Cayley graph.  Recall that the pair $(G,S)$ is called \textbf{almost convex}$(k)$ (denoted AC$(k)$) if there exists an integer $N(k)>0$ such that for every $n$ and for every pair of elements $g_1,g_2\in \Sigma(n)$, if $d_{S}(g_1,g_2)\leq k$, then there exists a path between $g_1$ and $g_2$ in $\beta(n)$ of length less than $N(k)$.  We say that $(G,S)$ is \textbf{almost convex (AC)} if it is AC$(k)$ for every $k$, and recall that Cannon proved that  AC$(2)$ implies AC. For the proof of this fact, and of other general properties of almost convexity, we refer the reader to \cite{CAC87}.
\section{2-Step Solvable non-nilpotent Lie Groups}\label{Solv}

The manifolds $X_A$ introduced in Section \ref{tbundles} are naturally compact quotients of certain two-step solvable Lie groups.  Given an element $A\in GL(n,\R)$, if $A$ is in the image of the exponential map $\exp:\mathfrak{gl}_n\rightarrow GL(n,\R)$, we can form the Lie group $\Gamma_A=\R^n\rtimes_{A^t}\R$ where the $\R$-factor on the right acts on $\R^n$ as multiplication by $A^t$.  If $A \in SL(n,\Z)$ and $A$ lies on a one-parameter subgroup as above, then the group $G_A=\Z^n\rtimes_A\Z$ sits naturally inside of $\Gamma_A$ as a cocompact lattice, namely as the subgroup of integer points in $\R^n$ and $\R$.
\subsection{ Jordan Blocks and One-parameter Subgroups}

If $\lambda \in \R_+$, we denote by $B(n, \lambda)$ the Jordan block of dimension $n$ with eigenvalue $\lambda$.  If $\lambda\in \C\setminus \R_+$, denote by $B(2n,\lambda, \overl{\lambda})$ the \emph{real} Jordan block of dimension $2n$ with eigenvalues $\lambda,\overl{\lambda}$. Concretely, we have \[B(n,\lambda)=\left (\begin{array}{cccc}\lambda& 1&\cdots&0\\ 
0&\lambda & \ddots & \vdots\\ 
\vdots& \vdots &\ddots&1\\  
0 &\cdots&0&\lambda

\end{array}\right)
\] for $\lambda \in \R$ and \[B(2n,\lambda, \overl{\lambda}) = \left (\begin{array}{c|c|c|c}P& I_2&\cdots&0\\ \hline
0&P & \ddots & \vdots\\ \hline
\vdots& \vdots &\ddots&I_2\\  \hline
0 &\cdots&0&P

\end{array}\right)\] for $\lambda \in \C\setminus \R_+$, where $I_2$ is the $2\times2$ identity matrix and if $\lambda=a+bi\mbox{, with $a,b\in \R$}$, then $P$ is the $2\times2$ matrix of the form\[P=\left(\begin{array}{cc} a &-b\\b&a\end{array}\right).\]

It is well-known that a matrix $A\in GL(n,\R)$ lies on a one-parameter subgroup if and only if in the Jordan block decomposition of $A$, blocks corresponding to negative eigenvalues come in pairs.  Given a matrix $A$, define $J(A)$ to be the \textbf{absolute Jordan form} of $A$ (Cf. \cite{FarbMosh00}), \emph{i.e.} the matrix obtained from the Jordan decomposition of $A$ by replacing the diagonal entries of negative eigenvalue blocks which cannot be paired by their absolute values:

\[\left (\begin{array}{cccc}\lambda& 1&\cdots&0\\ 
0&\lambda & \ddots & \vdots\\ 
\vdots& \vdots &\ddots&1\\  
0 &\cdots&0&\lambda

\end{array}\right)\rightarrow
\left (\begin{array}{cccc}|\lambda|& 1&\cdots&0\\ 
0&|\lambda| & \ddots & \vdots\\ 
\vdots& \vdots &\ddots&1\\  
0 &\cdots&0&|\lambda|

\end{array}\right)
\]
Then $J(A)$ always lies on a one-parameter subgroup, and we may form the Lie group $\Gamma_{J(A)}=\R^n\rtimes_{J(A)^t}\R$. For $\lambda \in \R_+$, the one-parameter subgroup containing a Jordan block $B=B(n,\lambda)$ has the following form
\[B(t)= \left (\begin{array}{cccc}\lambda^t& \lambda^{t-1}q_2(t)&\cdots&\frac{\lambda^{t-n+1}}{(n-1)!}q_n(t)\\ 
0&\lambda^t & \ddots & \vdots\\ 
\vdots& \vdots &\ddots&\lambda^{t-1}q_2(t)\\  
0 &\cdots&0&\lambda^t

\end{array}\right)\]where $q_2(t)=t$, and in general for $k>2$, $q_k(t)$ is a monic polynomial of degree $k-1$ which vanishes at $t=1$.  The corresponding one-parameter subgroup for a Jordan block $B=B(2n, \lambda,\overl{\lambda})$ has the form \[B(t) = \left (\begin{array}{c|c|c|c}\|\lambda\|^t\cdot P(t)& \|\lambda\|^{t-1}q_2(t)\cdot I_2&\cdots&\frac{\|\lambda\|^{t-n+1}}{(n-1)!}q_n(t)\cdot I_2\\ \hline
0&\|\lambda\|^t\cdot P(t)& \ddots & \vdots\\ \hline
\vdots& \vdots &\ddots&\|\lambda\|^{t-1}q_2(t)\cdot I_2\\  \hline
0 &\cdots&0&\|\lambda\|^t\cdot P(t)

\end{array}\right)
\] with the $q_k(t)$ the same as above. Here $P(t)$ a $2\times2$ matrix of the form \[P(t)=\left(\begin{array}{cc}\cos(\theta t)&-\sin(\theta t)\\
\sin(\theta t)&\cos(\theta t)\end{array}\right),\] where $\theta=\arg(\lambda)$.  If we have a pair of blocks with negative real eigenvalue $\lambda$, this belongs to a one-parameter subgroup of the latter type, with $\theta=\pi$. In each case, $B(1)=B$.  

With this notation, we can conjugate $A$ over $\R$ so that it has block diagonal form.  Then $J(A)$ has block diagonal form \[J(A)=\left(\begin{array}{cccc}
					B_1 &0&\cdots&0\\
					0& B_2&\cdots&0\\
					0&\cdots&\ddots&0\\
					0&\cdots&0&B_k	

				\end{array}\right)\] 
where each $B_i$ is a Jordan block as above.  $J(A)$ then lies on a one-parameter subgroup $J(A)^t$ of the form  \[J(A)^t=\left(\begin{array}{cccc}
					B_1(t) &0&\cdots&0\\
					0& B_2(t)&\cdots&0\\
					0&\cdots&\ddots&0\\
					0&\cdots&0&B_k(t)	

				\end{array}\right)\]
where each $B_i(t)$ is a corresponding one-parameter subgroup.   
\subsection{Isometries and Lattices of $\Gamma_{J(A)}$}
If we choose a left-invariant Riemannian metric on $\Gamma_{J(A)}$, then $\Gamma_{J(A)}$ itself can be embedded  in the identity component of the full isometry group $\text{Isom}(\Gamma_{J(A)})$ where it acts by left multiplication: every $g\in \Gamma_{J(A)}$ defines a map \begin{align*}l_g:\Gamma_{J(A)}&\rightarrow \Gamma_{J(A)}\\
h&\mapsto g\cdot h
\end{align*}
Using the standard coordinates on $\R^n\times \R$, define a left invariant metric on $\Gamma_{J(A)}$ by choosing the standard inner product at the identity.  One can check that the full isometry group with respect to this metric contains a subgroup isomorphic to $(\Z/2\Z)^k$ where $k$ is the number of Jordan blocks. The generator of each factor is the map $S_i:\R^n\rightarrow \R^n$ defined on basis elements $e_j$ by \[S_i(e_j)=\left\{\begin{array}{cl}
-e_j, &\mbox{ $e_j$ belongs to the $i$th Jordan block,}\\
e_j, & \mbox{ else}\end{array}\right.\]
We extend these to $\Gamma_{J(A)}$ by acting as the identity on the $\R$-factor.  It is clear that $S_i$ is an involution, and for all $g\in \Gamma_{J(A)}$ the derivative of $S_i$ commutes with the derivative of the left multiplication map $l_{g^{-1}}:\Gamma_{J(A)}\rightarrow \Gamma_{J(A)}$, hence it preserves the metric everywhere.  The importance of these extra elements of the isometry group is that even if $A$ does not lie on a one-parameter subgroup, the group $G_A\cong\pi_1(X_A)$ is a lattice in Isom$(\Gamma_{J(A)})$: 
\begin{lemma} \label{embed}$G_A\cong\pi_1(X_A)$ is a lattice in Isom$(\Gamma_{J(A)})$.

\end{lemma}
\footnotesize PROOF:
\normalsize

Recall the presentation for $G_A$ as an HNN-extension:
\[G_A = \langle a_1,\ldots,  a_n ,t_0 |\mbox{ }t_0  a_it_0^{-1}=Aa_i,\mbox{ $1\leq i \leq n$};\mbox{ } [a_i,a_j], \mbox{ $1\leq i,j\leq n$}\rangle.\]

For some product $S=S_{i_1}\cdots S_{i_n}$ we have that $A= SJ(A)$.  Let $e_i$ denote the $i$th standard basis vector of $\R^{n+1}$.  For $1\leq i\leq n$, set $v_i=l_{e_i}$ and define $t=S\cdot l_{e_{n+1}}$.  We consider the group $G = \langle v_1,\ldots, v_n, t\rangle\subseteq$ Isom$(\Gamma_{J(A)})$.  We will show that $G\cong G_A$.  First we check that the same relations in hold in $G$ that do in $G_A$.  Let $(\vec{x}, s)$ be an arbitrary element of $\Gamma_{J(A)}$.  Clearly we have \[(v_i v_j)\cdot(\vec{x},s) =(x+e_j+e_i,s)=(v_j v_i)\cdot(\vec{x},s)\]
hence all $v_i, v_j$ commute. For $1\leq i \leq n$ we also compute 
\begin{align*}
(tv_it^{-1})\cdot(\vec{x},s)&=(tv_i)\cdot(A^{-1}\vec{x},s-1)=t\cdot(A^{-1}\vec{x}+e_i,s-1)=(A\cdot A^{-1}\vec{x} +Ae_i,s-1+1)\\&= (\vec{x}+Ae_i,s)=(Av_i)\cdot (\vec{x},s).
\end{align*}

We can define a surjective map $\phi:G_A\rightarrow G$ sending $\phi: a_i\mapsto v_i$ and $\phi: t_0\mapsto t$. Note that from the action of $G$ on $\Gamma_A$, it is clear that $\langle v_1,\ldots, v_n\rangle$ is a free abelian subgroup of rank $n$, and that $\phi$ therefore maps $\langle a_1,\ldots, a_n\rangle$ isomorphically onto $\langle v_1,\ldots, v_n\rangle$. It remains to show that $\phi$ is injective on all of $G_A$. Suppose $\phi(g)=1_G$, and write $g=s_1\cdots s_k$ where the $s_i$ are generators of $G_A$.  Observe that the only generator of $G$ whose action on $\Gamma_A$ affects the last coordinate is $t$.  Since $\phi(g)$ takes $0\in \R^{n+1}$ to $0$, we must therefore have that the exponent sum of all $t$'s appearing in $\phi(g)$ is $0$.  This also implies that the exponent sum of $t_0$'s appearing in $g$ is 0, and hence $g$ is a product of conjugates of the $a_i$.  In particular, $g$ is an element of $\langle a_1,\ldots, a_n\rangle\cong \Z^n$.  By the remark above,  this implies $g=1_{G_A}$.  \qed

\bigskip

We now give the proof of Theorem \ref{main2}.

\bigskip
\footnotesize PROOF
\normalsize of Theorem \ref{main2}:

By Lemma \ref{embed}, we can think of $G=G_A$ as a cocompact lattice of Isom$(\Gamma_{J(A)})$.  We will assume that $A$ (no longer in $SL(n,\Z)$) is  a matrix in Jordan normal form.  $G$ is then identified with a discrete subgroup of $\Gamma_{J(A)}\rtimes(\Z/2)^k$ and we can write each element $g\in G$ uniquely as $g=d(g)s(g)$ with $d(g)\in \Gamma_{J(A)}$ and $s(g)\in (\Z/2)^k$. The translation component $d(g)$ will further be identified with the element of $\Gamma_{J(A)}$ that it represents.  In this way, $d(g)$ is a vector in $\R^{n+1}$, so we will speak of the coordinates of $d(g)$ as pertaining to particular eigenspaces.  The \emph{last} coordinate will always denote the $\R$-factor in the semidirect product defining $\Gamma_{J(A)}\cong\R^n\rtimes_{J(A)^t}\R$. 

Let $S=S^{-1}$ be some generating set for $G$.  Let $G_0$ be the subgroup $G\cap\mbox{Isom}_0(\Gamma_{J(A)})$.  Then if $s\in G\setminus G_0$ then $s^2\in G_0$.  We observe that if $s\in S$ is a generator with maximal last coordinate, then $s^2$ has length 2 in the word norm and has maximal last coordinate among all length 2 words.  Let $t_0$ be the maximal last coordinate occurring among all generators, and set $t^*=2t_0$.   

Let $w_0=(a_1,\ldots,a_n,t^{*})$ be the square of some generator with maximal last coordinate. Since $A$ has some eigenvalue not on the unit circle, it has one eigenvalue  with norm greater than one, and one with norm smaller than one.  Without loss of generality, we may decompose any $\textbf{v}\in \R^n$ as $\textbf{v}=(v_1,v_2,v_3)$ where $v_1$ belongs to the eigenspace $V_\lambda$ of a maximal-norm eigenvalue $\lambda$, $v_2$ belongs to the eigenspace $V_\mu$ of a minimal-norm eigenvalue $\mu$, and $v_3$ belongs to the complement of the first two.  Since $G \cap \R^n\times \{0\}$ is a rank $n$ lattice, we can find elements $\textbf{x}=(x_1,x_2,x_3,0)$ and $\textbf{y}=(y_1,y_2,y_3,0)$ in $G$ and we may assume that $x_1$ and $y_2$ are both positive vectors. 
 
 For each $k>0$, consider the following elements of $G$:
 \[X_k = w_0^k\cdot\textbf{x}\cdot w_0^{-k}\]
 \[Y_k = w_0^{-k}\cdot\textbf{y}\cdot w_0^{k }\]
 Then $X_k,Y_k$ are again elements of  $G \cap \R^n\times \{0\}$, and hence $X_k\cdot Y_k=Y_k\cdot X_k$.  Now for any $0\leq j \leq k$, we define two further elements:
 \[\sigma_j = X_k\cdot Y_k\cdot w_0^{-j}\]
 \[\tau_j = Y_k\cdot X_k \cdot w_0^{j}\]
 
By construction, $\sigma_j$ and $\tau_j$ are both within $2j$ of $X_k\cdot Y_k$, hence within $4j$ of eachother. We also have:\[\norm{\sigma_j}\mbox{,}\norm{\tau_j} \leq8k+ \norm{\textbf{x}}+\norm{\textbf{y}}-2j\]
If $k,j$ are sufficiently large then $2j$ dominates the two middle terms on the right and we can instead write 
\[\norm{\sigma_j}\mbox{,}\norm{\tau_j} \leq8k-2J\]
where $J$ is a constant which depends only on $\norm{\textbf{x}}$, $\norm{\textbf{y}}$ and $j$.  We will assume that $(G,S)$ is almost convex and show that for some fixed choice of $j$ to be specified later, $N(4j)$ cannot exist.   

Observe that the last coordinate of $\sigma_j$ is less than $0$ while the last coordinate of $\tau_j$ is greater than $0$.   Therefore any path between them in the Cayley graph for $G$ must pass through some group element whose last coordinate has absolute value less than or equal to $t^{*}$. In order to reach a contradiction, it suffices to show that any product of elements in $B(8k-2J)$ whose last coordinate is bounded by $t^{*}$ in absolute value can be made arbitrarily far from $X_k\cdot Y_k$ by choosing $k$ large enough.  Since $\sigma_k$ and $\tau_k$ are a bounded distance from $X_k\cdot Y_k$, this will prove the theorem.  

In order to prove the result, we first estimate the effect of taking products of group elements on a particular Jordan subspace of dimension $m$.   Let $v_1,\ldots, v_k$ be vectors in a Jordan subspace $V_\lambda$ associated to the eigenvalue $\lambda$.  If $\lambda$ is negative or complex we consider the subspace corresponding to $\lambda$ and $\overl{\lambda}$.  We may assume that this subspace corresponds to the first $m$ basis vectors $\{e_1,\ldots, e_m\}$ (resp. $2m$ basis vectors $\{e_1,e_1',\ldots, e_m,e_m'\}$).   As noted above, with respect to this basis the one parameter subgroup $B(t)$ is upper triangular (resp. block upper triangular).  We will be interested in the norm of the projection  of $V_\lambda$ onto the subspace $\langle e_m\rangle$ (resp. $\langle e_m,e_m'\rangle$). Define $\pi_m$ to be this projection operator.  

Consider \[w=(v_1,-,t_1)\cdots(v_k,-,t_k)= (x,-,t_1+\cdots +t_k),\] and let $T_j= \sum_{i=1}^j t_i$ be the $j^{th}$ partial sum. Then we have

\begin{lemma}\label{expgrow} For any such $w$, the following inequality holds \[(\ast)\norm{\pi_m(x)}\leq \norm{\pi_m(v_1)}+\|\lambda^{T_1}\|\cdot\norm{\pi_m(v_2)}+\cdots+\|\lambda^{T_k}\|\cdot\norm{\pi_m(v_k)}.\]
\end{lemma}

\footnotesize{PROOF}: 
\normalsize
 If $\lambda\in \R_+$, the Jordan block of type $B(m,\lambda)$ lies on a one parameter subgroup of the form:
\[B(t)= \left (\begin{array}{cccc}\lambda^t& \lambda^{t-1}q_2(t)&\cdots&\frac{\lambda^{t-m+1}}{(m-1)!}q_m(t)\\ 
0&\lambda^t & \ddots & \vdots\\ 
\vdots& \vdots &\ddots&\lambda^{t-1}q_2(t)\\  
0 &\cdots&0&\lambda^t

\end{array}\right).\]
Observe that $\pi_m(B(t)\cdot v)=\lambda^t\pi_m(v)$. Then  \begin{align*}\pi_m(x) &= \pi_m\left(v_1+B(T_1)\cdot v_2+\cdots +B(T_k)\cdot v_k\right)\\
&= \pi_m(v_1)+\pi_m(B(T_1)\cdot v_2)+\cdots +\pi_m(B(T_k)\cdot v_k)\\
&=\pi_m(v_1)+\lambda^{T_1}\pi_m(v_2)+\cdots \lambda^{T_k}\pi_m(v_k).\end{align*}
Taking norms of each side and applying the triangle inequality, we deduce $(*)$ in this case.  

If $\lambda \in \C\setminus \R_+$, then the Jordan block of type $B(2m, \lambda,\overl{\lambda})$ lies on a one-parameter subgroup of the form:
 \[B(t) = \left (\begin{array}{c|c|c|c}\|\lambda\|^t\cdot P(t)& \|\lambda\|^{t-1}q_2(t)\cdot I_2&\cdots&\frac{\|\lambda\|^{t-m+1}}{(m-1)!}q_m(t)\cdot I_2\\ \hline
0&\|\lambda\|^t\cdot P(t)& \ddots & \vdots\\ \hline
\vdots& \vdots &\ddots&\|\lambda\|^{t-1}q_2(t)\cdot I_2\\  \hline
0 &\cdots&0&\|\lambda\|^t\cdot P(t)

\end{array}\right)
\]where $P(t)$ is a $2\times 2$ rotation matrix of the form \[P(t)= \left(\begin{array}{cc}\cos(\theta t) &-\sin(\theta t)\\
\sin(\theta t) & \cos(\theta t)
\end{array}\right)\]  for $\theta =\arg(\lambda)$.  In this case, $\pi_m(B(t)\cdot v)=\|\lambda^t\|\cdot P(t)\cdot \pi_m(v)$
Since rotation matrices preserve norms we have
 \begin{align*}\norm{ \pi_m(x)} &= \norm{ \pi_m(v_1)+\pi_m(B(T_1)\cdot v_2)+\cdots +\pi_m(B(T_k)\cdot v_k)}\\
 &\leq \norm{\pi_m(v_1)}+\|\lambda^{T_1}\|\cdot \norm{P(T_1)\cdot \pi_m(v_2)}+\cdots \|\lambda^{T_k}\|\cdot\norm{P(T_k)\cdot \pi_m(v_k)}\\
 &\leq  \norm{\pi_m(v_1)}+\|\lambda^{T_1}\|\cdot\norm{ \pi_m(v_2)}+\cdots \|\lambda^{T_k}\|\cdot\norm{ \pi_m(v_k)}. \end{align*}

\QED

Set $N= 8k-2J$.  We will take a product of $N$ generators and assume that the last coordinate of the product has absolute value $\leq t^{*}$:\[\prod_{i=1}^{N}s_i = \prod_{i=1}^N(a_i,b_i,c_i,t_i)=(v_1,v_2,v_3,t)\]
where, again, we let $T_j = \sum_{i=1}^jt_i$ denote the $j^{th}$ partial sum. 

We have two cases to consider. In the first case, $T_j$ will be at most zero for more than half of the $j$, and in the second, we assume the opposite holds. We will prove that, in the first case, the product of $N$ generators has $\norm{\pi_m(v_1)}$ arbitrarily small as compared with the same component of $X_kY_k$.  In the second case, the same statement holds instead for $v_2$, and the proof is analogous so we omit it.  

Recall that $V_\lambda$ is the Jordan subspace containing $v_1$. Let $a^{*}$ be a vector occurring in the $V_\lambda$ component of some generator such that $\|\pi_m(a^*)\|$ is as large as possible.  Recall that as above, $t_0=t^{*}/2$ is the greatest value occurring in the $\R$-component for any generator.

Suppose $\dim(V_\lambda)=m$.  Since $T_j\leq0$ for more than half of the $j$, when $N$ is divisible by 4 the norm of $v_2$ is maximized when in the sequence of partial sums, each of $0, \pm t_0, \pm 2t_0,\ldots, \pm N/4-1$ occurs twice,  $\pm N/4$ once, and the $V_\lambda$ component of every factor is $a^{*}$.  Gathering each of these terms and counting 1 for all negative powers of $\lambda$, by $(*)$ we obtain:
\begin{align*}
\norm{\pi_m(v_1)} &\leq \norm{\pi_m(a^*)}\left\{2\sum_{j=0}^{N/4-1}(\|\lambda\|^{jt_0}+1)+\lambda^{(N/4)t_0}+1 \right\}\\
&\leq  \norm{\pi_m(a^*)}\left\{2\sum_{j=0}^{N/4}\|\lambda\|^{jt_0}+N/2\right\}\\
&\leq \norm{\pi_m(a^*)}\left\{C\|\lambda\|^{(N/4+1)t_0}+N/2\right\}
\end{align*}
where the last inequality was obtained by applying the binomial theorem, and $C=\frac{2}{\|\lambda\|^{t_0}-1}$. On the other hand we can estimate $X_k\cdot Y_k$ as follows:
 \[X_k = w_0^k\cdot(x_1,x_2,x_3,0)\cdot w_0^{-k} = (J(A)^{kt^*}\cdot x_1,J(A)^{kt^*}\cdot x_2, J(A)^{kt^*}\cdot x_3,0)\]
\[Y_k = w_0^{-k}\cdot(y_1,y_2,y_3,0)\cdot w_0^{k} = (J(A)^{-kt^*}\cdot y_1,J(A)^{kt^*}\cdot y_2, J(A)^{kt^*}\cdot y_3,0)\]
We can bound the norm of the projection of the $V_1$ component from below:
\[\norm{\pi_m\left(J(A)^{kt^*}\cdot x_1+J(A)^{-kt^*}\cdot y_1\right)}\geq\left| \|\lambda\|^{kt^*}\cdot \norm{\pi_m(x_1)}-\|\lambda\|^{-kt^*}\cdot\norm{\pi_m(y_1)}\right|\]

For $k$ sufficiently large we can clearly remove the outer absolute value bars on the right.  Putting this together with the bound for $\norm{\pi_m(v_1)}$ we obtain: 
\begin{align*}
\norm{\pi_m\left(J(A)^{kt^*}\cdot x_1+J(A)^{-kt^*}\cdot y_1\right)}-\norm{\pi_m(v_1)}&\geq \|\lambda\|^{kt^*}\cdot \norm{\pi_m(x_1)}-\|\lambda\|^{-kt^*}\cdot\norm{\pi_m(y_1)}\\&\hspace{1cm}- \norm{\pi_m(a^*)}\left\{C\|\lambda\|^{(N/4+1)t_0}+N/2\right\}\\
&\geq\|\lambda\|^{kt^*}\cdot\norm{\pi_m(x_1)}-C\|\lambda\|^{(N/4+1)t_0}\cdot\norm{\pi_m(a^*)} +O(1).
\end{align*}
Now recall that $t_0=t^{*}/2$, $N=8k-2J$ and therefore \[(N/4+1)t_0=\left(\frac{8k-2J}{4}+1\right)(t^{*}/2)=(k-J/4+1/2)t^*.\] From this we see that:
\[\norm{\pi_m\left(J(A)^{kt^*}\cdot x_1+J(A)^{-kt^*}\cdot y_1\right)}-\norm{\pi_m(v_1)}\geq\|\lambda\|^{kt^*}\left(\norm{\pi_m(x_1)}-C\|\lambda\|^{(-J/4+1/2)t^*}\cdot\norm{\pi_m(a^*)}\right) +O(1).\]

Since $J>0$, we may choose it in advance large enough so that \[\norm{\pi_m(x_1)}-C\|\lambda\|^{(-J/4+1/2)t^*}\cdot\norm{\pi_m(a^*)}>0.\]  Then as $k\rightarrow\infty$, we can make the difference arbitrarily large.  \QED

\bibliographystyle{plain}
\bibliography{Growth}

\end{document}